\documentclass[11pt,reqno, letterpaper]{amsart}
\usepackage{amsmath,amsfonts,amsthm,amscd,amssymb,graphicx,xcolor,tikz}
\usepackage{epstopdf}
\numberwithin{equation}{section}
\usepackage{fullpage}
\usepackage{color}
\usepackage[normalem]{ulem}

\newcommand\s{\sigma}

\renewcommand\a{\alpha}
\renewcommand\b{\beta}
\renewcommand\o{\omega}
\newcommand {\la}{\lambda}

\def\l{\lambda}
\def\eps{\varepsilon }
\def\e{\varepsilon}

\renewcommand\a{\alpha}

\newcommand\sig{\sigma}
\newcommand\si{\sigma}
\newcommand\pa{\partial}

\renewcommand\b{\beta}
\renewcommand\o{\omega}
\newcommand \lam {\lambda}
\newcommand\T{\mathbb T}
\newcommand\R{\mathbb R}

\def\eps{\varepsilon}
\def\e{\varepsilon}

\def\l{\lambda}

\newcommand\br{\begin{remark}}
\newcommand\er{\end{remark}}
\newcommand\brs{\begin{remarks}}
\newcommand\ers{\end{remarks}}
\newcommand\bp{\begin{pmatrix}}
\newcommand\ep{\end{pmatrix}}
\newcommand{\be}{\begin{equation}}
\newcommand{\ee}{\end{equation}}
\newcommand{\ben}{\begin{equation*}}
\newcommand{\een}{\end{equation*}}
\newcommand\ba{\begin{equation}\begin{aligned}}
\newcommand\ea{\end{aligned}\end{equation}}

\newcommand{\bap}{\begin{app}}
\newcommand{\eap}{\end{app}}
\newcommand{\begs}{\begin{exams}}
\newcommand{\eegs}{\end{exams}}
\newcommand{\beg}{\begin{example}}
\newcommand{\eeg}{\end{exaplem}}
\newcommand{\bpr}{\begin{proposition}}
\newcommand{\epr}{\end{proposition}}
\newcommand{\bt}{\begin{theorem}}
\newcommand{\et}{\end{theorem}}
\newcommand{\bc}{\begin{corollary}}
\newcommand{\ec}{\end{corollary}}
\newcommand{\bl}{\begin{lemma}}
\newcommand{\el}{\end{lemma}}
\newcommand{\bd}{\begin{definition}}
\newcommand{\ed}{\end{definition}}

\newcommand{\CalF}{\mathcal{F}}

\newcommand{\calL}{\mathcal{L}}

\newcommand{\calR}{\mathcal R}

\newcommand{\RR}{{\mathbb R}}

\newcommand{\ZZ}{{\mathbb Z}}
\newcommand{\TT}{{\mathbb T}}
\newcommand{\CC}{{\mathbb C}}

\newcommand {\ddt} {\frac{d}{dt}}
\newcommand{\va}{V^{\textup{app}}}
 \newcommand{\rap}{\mathcal{R}^{\textup{app}}}
\newcommand{\tu}{\tilde u}
\newtheorem{theorem}{Theorem}[section]
\newtheorem{proposition}[theorem]{Proposition}
\newtheorem{corollary}[theorem]{Corollary}
\newtheorem{lemma}[theorem]{Lemma}

\theoremstyle{remark}
\newtheorem{remark}[theorem]{Remark}
\newtheorem{remarks}[theorem]{Remarks}
\theoremstyle{definition}
\newtheorem{definition}[theorem]{Definition}

\newtheorem{example}[theorem]{Example}

\newcommand\cG{{\mathcal  G}}

\newcommand\cL{{\mathcal  L}}

\newcommand\cF{{\mathcal  F}}

\newcommand{\bbr}{\mathbb{R}}

\newcommand{\bbR}{{\mathbb{R}}}

\newcommand{\beq}{\begin{equation}}
\newcommand{\eeq}{\end{equation}}

\title{Nonlinear instability of rolls in the 2-dimensional generalized Swift-Hohenberg equation}

\date{}

\author{Myeongju Chae}
\address{Department of Applied Mathematics,  Hankyong National University, Anseong-si, Gyeonggi-do,
Korea}
\email{mchae@hknu.ac.kr}
\thanks{Research of M. C was supported by  the National Research Foundation of Korea (NRF) grant funded by the Korea government (No. RS-2023-00279920)}
\author{Soyeun Jung}
\address{Division of International Studies, Kongju National University, Gongju-si, Chungcheongnamdo, Korea}
\email{soyjung@kongju.ac.kr}
\thanks{Research of S.J. was supported by the National Research Foundation of Korea (NRF) grant funded by the Korea government (No. 2022R1F1A1074414). }

\begin{document}

\begin{abstract} Within the framework developed in \cite{Gr, JLL, RT1}, we rigorously establish the nonlinear instability of roll solutions to the two-dimensional generalized Swift-Hohenberg equation (gSHE). Our analysis is based on spectral information near the maximally unstable Bloch mode, combined with precise semigroup estimates.  We construct a certain class of small initial perturbations that grow in time and cause the solution to deviate from the underlying roll solution within a finite time. This result provides a clear  transition from spectral to nonlinear instability in a genuinely two-dimensional setting, where the Bloch parameter $\sigma$ ranges over an unbounded domain. 
\end{abstract}

\date{}
\maketitle


\section{Introduction}

This paper investigates the nonlinear instability of roll solutions to the two--dimensional gSHE with quadratic--cubic nonlinearity:
\be \label{2gSHE}
\partial_t u=-(1+\partial_{x}^2+\partial_y^2)^2u+\varepsilon^2u+bu^2-su^3, \qquad u\in \R. 
\ee
Here, $\e \in \R$ is a control parameter, and $b, s \in \RR$ are nonzero constants satisfying $b \neq 0$ and $27s-38b^2>0$. 

The Swift–Hohenberg equation (SHE), in its classical form with $b=0$, was  originally introduced as a model for the onset of spatially periodic structures near pattern-forming instabilities, especially in Rayleigh–Bénard convection (\cite{SH}). Since then, it has served as a canonical fourth-order partial differential equation for studying pattern formation in various physical, chemical, and biological systems (\cite{CH}). Its mathematical simplicity and rich dynamics have made it a prototypical model for analyzing bifurcation, stability, and dynamics of spatially periodic patterns. In two or more spatial dimensions, one fundamental class of solutions is given by roll solutions, one-dimensional periodic structures extended in the transverse direction, which are widely regarded as a basic example in the study of pattern formation (\cite{H}).

To capture more realistic nonlinear effects, several generalizations of the SHE have been proposed. In particular, the gSHE, which includes both quadratic and cubic nonlinear terms as in \eqref{2gSHE}, introduces asymmetric nonlinear effects and supports a wider range of pattern-forming behaviors than the classical SHE; e.g., see \cite{ALBKS, BK, CG, GX, HMBD, K, R} for detailed accounts of pattern diversity in the gSHE. Among these, spatially periodic patterns remain a central object of study, not only because of their role in pattern formation theory, but also due to because their stability and instability exhibit additional complexity when the reflection symmetry ($u \rightarrow -u$) is broken.


The spectral and nonlinear stability of spatially periodic solutions has been extensively studied in the SHE. The fundamental works \cite{M1, M2} rigorously characterized the spectrum of the linearized operator around roll solutions in two dimensions and analyzed their spectral stability and instability;  in particular, the stability boundaries are determined by the classical Eckhaus and zigzag instabilities. The nonlinear stability of spatially periodic solutions within the Eckhaus stable regime was first investigated in \cite{S1, S2} for one dimension and later extended to two dimensions in \cite{U}.

Similar studies have also been carried out in other models that exhibit spatially periodic patterns. For example, in the Brusselator model (\cite{SZJV}), the spectral stability and instability of spatially periodic  wave have been investigated using techniques from \cite{M2}. In addition, under the assumption of spectral stability, a general framework for deriving both linear and nonlinear stability has been developed and applied to various classes of evolution equations (e.g., see  \cite{JNRZ1, JZ} and the references therein).  These contributions form a well-developed body of theory on the stability of spatially periodic waves.

In contrast, nonlinear instability, that is, the nonlinear dynamics near a spectrally unstable spatially periodic pattern, remains less well understood. In particular, the transition from spectral to nonlinear instability poses significant analytical challenges, especially in multi-dimensional settings where mode dispersion and interactions are more intricate.

The purpose of this paper is to address this gap by establishing the nonlinear instability of roll solutions in the two-dimensional gSHE. Building on our previous spectral results in \cite{CJ}, we construct suitable perturbations and show their nonlinear growth, thereby demonstrating that solutions deviate from the underlying roll solution within a finite time.

Our analysis is based on a well-established strategy initiated in \cite{Gr}, where the growth properties of the linearized operator and the existence of an exponentially growing mode are combined to construct a nonlinear approximate solution that exhibits instability of the original nonlinear system. This approach was first implemented for two--dimensional Euler and Prandtl equations. Similar ideas also appear in \cite{RT1} and \cite{RT2}, which  establish the transversal instability of the line solitary water waves and the KdV solitary wave in two dimensions. See also \cite{JLL}, where nonlinear modulational instability of periodic traveling waves for several dispersive PDEs  is proved in one dimension, under both periodic and localized perturbations. 

Our goal in this introduction is to briefly summarize our previous results and to state the main theorem of the present paper. In \cite{J}, we constructed roll solutions for the one-dimensional gSHE and carried out their spectral analysis. In \cite{CJ}, we extended this to the two-dimensional setting and proved that all roll solutions are spectrally unstable throughout their existence regime. These results provide the spectral groundwork necessary for the nonlinear instability established in this paper. 

\smallskip


\textbf{Notation.} We introduce some basic notation used throughout the paper. 
For a nonnegative integer $p$, let $H^p(\RR^2)$ be the Sobolev space 
of order $p$ on $\RR^2$ with norm $\|\cdot\|_{H^p(\RR^2)}$; 
in particular, $H^0 = L^2$. 
Define $\TT = [0,2\pi]$ and let $H^p(\TT)$ be $H^p([0,2\pi])$ with periodic 
boundary conditions $\partial_x^j u(0) = \partial_x^j u(2\pi)$ 
for $j=0,\dots,p-1$. 
Unless otherwise specified, $H^p(\RR^2)$ and $H^p(\TT)$ are taken in 
$(x,y)\in\RR^2$ and $x\in\TT$, respectively.  

Throughout the paper, $x$ and $y$ are spatial variables, $t$ is the temporal variable, and there are four parameters $\e$, $\o$, $\s_1$ and $\s_2$, in particular, $\s:=(\s_1, \s_2)$ is a Bloch wave vector.
For $A,B>0$, we write $A\simeq B$ if there exists a uniform constant $C\ge1$ such that 
\[
  \tfrac1C A \le B \le C A.
\]
The symbol $C$ also denotes a generic uniform constant, possibly changing from 
line to line; dependence on parameters is made explicit when relevant.  The spectrum of an operator is denoted by $Spec$.

\subsection{Roll solutions}
A roll solution, as considered in this work, refers to a spatially periodic solution   of the form $$u_{roll}(t,x,y)=\tilde u(x),$$ which is independent of both $y$ and $t$. Such solutions bifurcate from the uniform state $u \equiv 0$, which becomes unstable as the control parameter $\e$ crosses a critical threshold (namely, $\e=0$). To understand the onset of instability, we linearize \eqref{2gSHE} about $u \equiv 0$, yielding 
$$\partial_t u=-(1+\partial_{x}^2+\partial_y^2)^2u+\varepsilon^2u.$$
Since we focus on solutions that are independent of $y$, we restrict our attention to perturbations of the form $u(t,x,y)=e^{\l t}e^{ikx}$, where $k\in \RR$ denotes the one-dimensional spatial wave number. Substituting into the linearized equation yields the dispersion relation $$ \l(k)=-(1-k^2)^2+\e^2.$$ This relation shows that the mode with wave number $k =\pm 1$ becomes unstable when $\e>0$, indicating a Turing-type bifurcation triggered by a finite-wavenumber instability. 

To simplify the analysis, we rescale the spatial variable by setting $x \rightarrow x:=kx$. With this change of variables, the equation \eqref{2gSHE} reads 
\be \label{2DgSH_intro}
\partial_t u=-(1+k^2\partial_{x}^2+ \partial_y^2)^2u+\varepsilon^2u+bu^2-su^3, \qquad u \in \RR^1.
\ee
The rigorous construction of such roll solutions with period $2\pi$ was carried out in our previous work \cite{CJ, J}. They have small amplitude, with leading-order behavior of size $\mathcal{O}(\e)$, as shown in the following theorem.   
\begin{theorem}[Existence of the rolls, \cite{CJ, J}] \label{existence_roll} Assume $27s-38b^2>0$ and let $\o$ be a parameter satisfying
\be \label{omega_intro}
k^2-1=2\e\o.
\ee
Then there exists an $\e_0>0$ such that for all $\e \in (0, \e_0]$ and all $\o \in [-\frac{1}{2}, \frac{1}{2}]$ there is a unique (up to translation) stationary $2\pi$--periodic roll solution $$u_{roll}(t, x,y)=\tilde u_{\e, \o}(x)$$ of \eqref{2DgSH_intro}, which is even in $x$ and bifurcating from the uniform state $u\equiv 0$.  These periodic solutions have the following expansion
\be \label{u_tilde}
\begin{split}
\tilde u_{\varepsilon, \omega}(x)
& = \e \frac{6 \sqrt{1-4\omega^2} }{ \sqrt{27 s-38 b^2}}\cos x \\
& \quad + \e^2 \Big[ \frac{18b (1-4\omega^2)}{27s-38b^2} -\frac{32b^2  \omega \sqrt{1-4\omega^2}}{ (27 s-38 b^2)\sqrt{27 s-38 b^2}} \cos x   +\frac{2b (1-4\o^2) }{27s-38b^2}\cos 2x  \Big]  \\
&  \quad + \mathcal{O}(\varepsilon^3).
\end{split}
\ee
In particular, if $\o=\pm \frac{1}{2}$ then $\tilde u_{\e, \o}(x) \equiv 0$.
\end{theorem}

Based on the parameter relation given in \eqref{omega_intro}, we will consider the equation in the following form throughout this paper: 
\be \label{2DgSH_w}
\partial_t u=-(1+(1+2\e \o)\partial_{x}^2+ \partial_y^2)^2u+\varepsilon^2u+bu^2-su^3, \qquad u \in \RR^1.
\ee



\subsection{Spectral instability} \label{Spectral instability} We now briefly review the spectral analysis of the roll solutions \eqref{u_tilde} carried out in \cite{CJ}. To this end, we first recall the Bloch--wave decomposition, which is a central tool not only for studying the linearized spectrum but also for deriving the semigroup estimates in Chapter \ref{Semigroup estimates}. 

Linearizing \eqref{2DgSH_w} about $\tilde u_{\varepsilon, \omega}$ gives the linear operator
\be \label{linearization about rolls}
\mathcal{L}_{\e, \o} :=-(1+(1+2\e \o)\partial_{x}^2+\partial_y^2)^2  +\e^2 +2b \tilde u_{\eps,\omega} -3s\tilde u_{\eps,\omega}^2
\ee
acting on $L^2 (\RR^2)$ with densely defined domain $H^4(\RR^2)$. Since $\tilde u_{\e, \o}$ is $2\pi$--periodic in $x$ and independent of $y$, so are the coefficients of $\mathcal{L}_{\e, \o}$. Then, by the standard Floquet--Bloch theory(\cite{JmZk, KP, M2}), any bounded eigenfunction $v(x, y)$ of $\mathcal{L}_{\e, \o}$ takes the form
\be \label{form of v}
v(x, y)=e^{i(\sigma_1 x+\sigma_2 y)} W(x, \s, \l), \quad W(x+2 \pi, \s, \l)=W(x, \s, \l),
\ee
where $\s:=(\s_1, \s_2) \in \RR^2$ are referred to as Bloch wave vectors throughout this paper. Here $W$ depends only on $x$, so it cannot belong to $L^2(\RR^2)$. As a result, the linear operator $\mathcal{L}_{\e, \o}$ admits no $L^2(\RR^2)$--eigenfunctions, and its entire spectrum is essential (see \cite[Chapter 3.3]{KP}). To analyze this essential spectrum, we introduce the family of Bloch operators below. 

By substituting \eqref{form of v} into \eqref{linearization about rolls}, we define the Bloch operator $B(\e, \o, \s)$ for each $\s=(\s_1, \s_2)$ by 
\begin{align}\label{Bloch_r}
\begin{aligned}
B(\eps, \omega, \s)
:& =e^{-i(\s_1x +\s_2 y)} \calL_{\e, \o}e^{i(\s_1x +\s_2 y)}\\
& =-\big[ 1+(1+2\e \o)(\partial_x+i\sigma_1)^2-\sigma_2^2 \big]^2  + \e^2 +2b \tilde u_{\eps,\omega} -3s\tilde u_{\eps,\omega}^2
\end{aligned}
\end{align}
acting on $L^2(\mathbb T)$ with a densely defined domain $H^4(\mathbb T)$. For each $\s \in \RR^2$, $L^2(\mathbb T)$--spectrum of $B(\e, \o, \s)$ consists entirely of discrete eigenvalues; that is, it is purely  point spectrum. Moreover, the operator $B(\e, \o, \s)$ is self--adjoint, and hence all its eigenvalues are real. Of particular importance is the fact that the full $L^2(\RR^2)$--spectrum of $\mathcal{L}_{\e, \o}$ can be described as the closure of the union of these  point spectra over all $\sigma \in \RR^2$:
\be \label{spectral identity}
Spec_{L^2(\RR^2)}(\mathcal{L}_{\e, \o})=closure \Big( \bigcup_{\s \in \RR^2} Spec_{L^2(\mathbb T)}(B(\e, \o, \s)) \Big). 
\ee
The interested reader can consult \cite[Section 2]{M2} for the proof of this spectral identity.

In \cite{CJ},  we identified a region in the Bloch parameter space $\s$ where instability was expected (see \cite[Section 3.1]{CJ} or the brief review in Section \ref{Linear instability of rolls}). Within this region, we studied the eigenvalue problem for the Bloch operator $B(\e, \o, \s)$ and showed that there exists a specific Bloch wave vector $\s^*$ such that, for sufficiently small $\e>0$ and all $\o \in [-\frac{1}{2}, \frac{1}{2}]$, the operator $B(\e, \o, \s)$ has a positive eigenvalue for every $\s$ sufficiently close to $\s^*$. This result is stated in the following theorem. Here, $\mathcal{S}_{\e, \o}$ denotes the set of all unstable Bloch wave vectors, i.e., those $\s \in \RR^2$ for which $B(\e, \o, \s)$ has a positive eigenvalue. 

\begin{theorem} [Spectral instability, \cite{CJ}]\label{spectral instability}  Assume $27s-38b^2>0$ and $b \neq 0$.  Let $\e_0$ be taken from the existence result of $\tilde u_{\e, \o}$ in Theorem \ref{existence_roll}. Then there exists  ${\tilde \e}_0 \in (0, \e_0]$ such that $\tilde u_{\e, \o}$ is spectrally unstable for all $\e \in (0, {\tilde \e}_0]$ and all $\o \in [-\frac{1}{2}, \frac{1}{2}]$. In particular, for any fixed $r>1$,
\begin{center}
$\{\s \in \RR^2~|~|\s-\sigma^*|\leq \mathcal{O}(\e^r) \} \subset \mathcal{S}_{\e, \o}$ \quad for all $\o \in [-\frac{1}{2}, \frac{1}{2}]$,
\end{center}
where $\sigma^*:=(-\frac{1}{2}, \frac{\sqrt{3}}{2})$.
\end{theorem}

This result contrasts with the classical SHE with $b=0$, where the instability of roll solutions is characterized by well-known instability boundaries in the parameter space (see \cite{M2} and the references therein). These boundaries, called the Eckhaus and zigzag instabilities, describe the regions in the plane of control parameters and wave numbers where roll patterns lose stability. More precisely, the Eckhaus instability occurs when perturbations act along the roll direction (\cite{E, TB}), while the zigzag instability arises when perturbations act in the transverse direction to the rolls (\cite{C, CH}).

In the two-dimensional gSHE, however, the inclusion of a quadratic nonlinearity ($b\neq 0$) breaks the reflection symmetry and leads to a different stability scenario. While in one dimension spatially periodic solutions satisfy the Eckhaus stability boundary established in \cite{J}, in two dimensions all roll solutions are spectrally unstable throughout the entire parameter regime.  Thus, unlike the classical SHE, the gSHE does not admit a stable region bounded by Eckhaus or zigzag curves. This shows that the quadratic term induces instability mechanisms that are genuinely two-dimensional and cannot be explained by one-dimensional analysis.



\subsection{Main result and strategy}  We now state our main theorem on the nonlinear instability of the roll solutions and outline the strategy of the proof. The theorem shows that  arbitrarily small perturbations can lead to a significant deviation from the roll solution within a finite time.



\begin{theorem} [Nonlinear instability] \label{nonlinear instability} 
 Let $s$, $b$ and $\tilde \e_0$ be chosen as in Theorem \ref{spectral instability}. Then there exists $\bar \e_0 \in (0, \tilde \e_0]$ such that for every $0 < \e< \bar \e_0$ and every $\o \in (-\frac{1}{2}, \frac{1}{2})$, the following nonlinear instability result holds.

There exists $\theta>0$ such that for every $\delta>0$ we can find an initial data $V_{\delta}(0, x, y)$
and a time  $T^{\delta}\sim |\ln \delta|$ such that
\be \label{goal_delta}
\|V_{\delta}(0, x, y)-\tilde u_{\e, \o}(x)\|_{H^4(\RR^2)} \leq \delta, 
\ee
and there exists a solution $V_{\delta}(t,x,y) \in H^2(\RR^2)$ of \eqref{2DgSH_w}
with initial data $V_{\delta}(0, x, y)$, defined on $[0, T^{\delta}]$, satisfying
\be \label{goal_theta}
\|V_{\delta}(T^{\delta}, x, y)- \tilde u_{\e, \o}(x)\|_{L^2(\RR^2)} \geq \theta.
\ee
\end{theorem}
\begin{proof} Given in Section \ref{Nonlinear instability}. 
\end{proof}

We note that the same type of nonlinear instability also holds in the spectrally unstable regime of the one-dimensional gSHE, namely for $\omega \in \big(-\tfrac{1}{2}, -\tfrac{1}{2\sqrt{3}}\big) \cup \big(\tfrac{1}{2\sqrt{3}}, \tfrac{1}{2}\big)$ (see \cite{J}).   Our contribution here is to establish this result in the genuinely two-dimensional setting, where the analysis is more delicate.

To analyze the nonlinear instability of roll solutions, it is essential to construct an exponentially growing solution to the linearized equation 
\be \label{LE}
V_t = \cL_{\e,\o}V,
\ee
where $\cL_{\e,\o}$ is defined in \eqref{linearization about rolls}. In our previous work \cite[Proposition 3.1]{CJ}, it was shown that certain spatially localized initial data $U$ lead to solutions satisfying the lower bound
\be
\| (e^{\cL_{\e, \omega} t } U)(x, y) \|_{H^2(\RR^2)} \gtrsim \frac{e^{\l(\s^*)t}}{(1+ t)^{\frac 1l}}
\notag
\ee
for some positive constant $l$, where $\l(\s^*)$ is an unstable eigenvalue associated with the mode $\s^*$ identified in Theorem \ref{spectral instability}. This demonstrates the linear instability of the roll solution. 

For nonlinear instability, however, such a lower bound is insufficient. To control the nonlinear terms, we need a more precise description of the dominant growth, including an upper bound on the linear evolution. In particular, we require the refined exponential estimate
\be \label{linearintro}
\| (e^{\cL_{\e, \omega} t } U)(x, y) \|_{H^2(\RR^2)} \simeq e^{\l(\bar \s) t} e(t), 
\ee
where $\l(\bar \s)$ is the maximal growth rate and $e(t)$ is a time-dependent prefactor satisfying $C_1(1+t)^{-1/l} \leq e(t) \leq C_2$ for some constants  $C_1$, $C_2>0$  (see Proposition \ref{liin} for the precise form of $e(t)$). This estimate captures the leading linear growth while keeping the nonlinear remainder under control. 

In addition, as indicated in \eqref{linearintro}, we need to focus on the most unstable mode $\bar \s$, the Bloch wave vector that maximizes the unstable eigenvalue of the Bloch operators. While in \cite{CJ} it was enough to construct growing solutions from an arbitrary unstable mode, here we must work with the maximizer $\bar \s$ in order to control the nonlinear terms in the perturbation analysis. To this end, in Section \ref{Linear instability of rolls} we revisit the linear instability analysis based on the maximizer $\bar \s$, rather than the previously used mode $\sigma^*$. 

In Section \ref{Semigroup estimates}, we establish semigroup estimates for the linearized evolution operator $e^{t\mathcal{L}_{\varepsilon,\omega}}$, which play a central role in the nonlinear analysis. These estimates show that the growth of the linear evolution is governed by the maximal unstable eigenvalue and are further used to control higher-order nonlinear terms through inhomogeneous bounds.

Lastly, in Section \ref{Nonlinear instability}, we prove the main theorem by constructing a solution that is decomposed into an approximate solution and a nonlinear remainder. Following the strategy of \cite{Gr, JLL, RT1}, the approximate solution is designed to match the true dynamics up to a certain order, while the remainder is controlled by energy estimates so that nonlinear effects remain subordinate to the leading exponential growth.

\subsection{Discussion}


It is worth noting two representative studies, Jin–Liao–Lin~\cite{JLL} and Rousset–Tzvetkov \cite{RT1}, which both employ the Grenier framework~\cite{Gr} to promote spectral instability to nonlinear instability.

The work \cite{JLL}, which studies the nonlinear modulational instability of periodic traveling waves in a one-dimensional dispersive model, shares structural features with the present paper. In both works, the linearized operators possess purely essential spectrum, but our two-dimensional setting gives rise to additional difficulties that do not appear in \cite{JLL}.

In one dimension the Bloch parameter $\s$ varies over a compact set, whereas in two dimensions it extends over an unbounded domain. Consequently, in our framework, the existence of a maximizer $\bar\sigma$ attaining the maximal growth rate $\lambda(\bar\sigma)$ is not immediate. We clarify this issue by showing that, outside a compact region in the $\s$-space, the spectrum of the Bloch operator $B(\varepsilon,\omega,\sigma)$ is uniformly bounded above by a negative constant (see Lemma \ref{sqrte}).

Another issue specific to the two-dimensional setting arises in establishing upper and lower bounds for the linear solution associated with the most unstable mode $\bar\sigma$ (see Proposition \ref{liin}). To obtain these bounds, a convenient approach is to employ a Taylor-type expansion of the unstable eigenvalue $\lambda(\sigma)$ near $\bar\sigma$, which in turn requires a certain regularity of $\lambda(\sigma)$ with respect to $\s$. In the spectral theory of periodic Schrödinger operators, it is well known that $\l(\s)$ is analytic in one dimension, whereas in higher dimensions it is in general only continuous (see~\cite[Theorem 5.5]{Ku} and~\cite{Wi}).
This lack of analyticity usually prevents the use of Taylor expansions in multidimensional settings. In our case, however, the local analyticity of $\lambda(\sigma)$ near any unstable mode is ensured by the Lyapunov–Schmidt reduction developed in \cite{CJ} (see Lemma \ref{analytic lambda_3}).

The work of Rousset--Tzvetkov~\cite{RT1} establishes the linear and nonlinear instability of line solitary water waves in two dimensions with respect to transverse perturbations. More precisely, they prove that one-dimensional solitary waves, which are (orbitally) stable under one-dimensional perturbations (Mielke \cite{M3}), become nonlinearly (orbitally) unstable when subjected to perturbations depending on the transverse variable. 
In their framework, the spectral problem is reduced to one-dimensional linearized operators $JL(k)$ parameterized by the transverse Fourier mode $k$, analogous to our $\sigma_2$, and the monotonic dependence of the Dirichlet-Neumann operator $G_{\e, k}$ on $|k|$ is crucial in proving the existence of the most unstable mode. Their analysis also involves a more intricate Hamiltonian system with several delicate analytical components, such as establishing long-time existence for the remainder equation.

In contrast to \cite{JLL, RT1, RT2}, our setting is dissipative rather than Hamiltonian. This dissipative structure considerably simplifies the analytical framework, such as the semigroup estimates (Section \ref{Semigroup estimates}) and the control of the remainder terms (Section \ref{Nonlinear instability}), and thus allows the Grenier framework to appear in transparent form. Nevertheless, to the best of our knowledge, our work provides the first application of the Grenier scheme to a dissipative equation. This is possible because in our case the spectral information of the Bloch operators, including the existence of the most unstable mode, is fully understood through Bloch-Floquet theory and Lyapunov-Schmidt reduction. Although the existence of an unstable mode implies linear instability, identifying the most unstable mode that derive the nonlinear instability is far from straightforward in general. In the context of dissipative systems, we also note that M.~Colombo et~al.~\cite{CO} recently established the linear instability of shear flows for the two-dimensional Navier–Stokes equations on a rectangular torus. In the presence of physical boundaries with no-slip conditions, the works \cite{BG}, \cite{GGN1}, \cite{GGN2}  have investigated spectral instability mechanisms arising from boundary layers. For further references on this topic, see the bibliography in \cite{CO}.

A possible direction for future work is to extend the nonlinear instability result to a broader class of initial perturbations. Indeed, Theorem \ref{nonlinear instability} is proved for a class of initial data designed to activate the most unstable part of the linearized operator. This guarantees maximal growth but relies on the specific structure of the perturbation. For comparison, the work of \cite{CS} established nonlinear instability for a wider class of perturbations in one-dimensional reaction–diffusion equations with periodic coefficients, under suitable spectral and nonlinear assumptions. Extending our approach to such more general perturbations, as well as to higher-dimensional settings, remains an interesting direction for future research.

\section{Linear instability of rolls} \label{Linear instability of rolls}

In this section, we construct an exponentially growing solution to the linearized equation \eqref{LE}. Our construction modifies the approach in \cite[Proposition 3.1]{CJ}, and provides an additional upper bound on the growing solution, which will be crucial for analyzing the nonlinear instability of the rolls solution. For this reason, we concentrate on the most unstable mode, namely the Bloch wave vector $\sigma$ corresponding to the largest unstable eigenvalue of the Bloch operator $B(\e, \o, \sigma)$. 
The existence of such a maximizer, together with the local real analyticity of the corresponding eigenvalue, is guaranteed by Lemma \ref{sqrte} and Lemma \ref{analytic lambda_3}.



Before presenting the lemmas, and to aid the reader’s understanding, we briefly recall the four regions in the Bloch parameter space $\s$ where unstable Bloch modes may occur (\cite[Section 3.1]{CJ}). Thanks to the spectral identify \eqref{spectral identity}, we can reduce Bloch wave vectors $\s \in \RR^2$ to ${\sigma} \in [-\tfrac{1}{2}, \tfrac{1}{2}] \times [0, \infty) $
 using that \( W \) is \( 2\pi \)-periodic in \( x \) and that \( B(\e, \o, \s) \) is even in \( \sigma_2 \). 
 Furthermore, we can further restrict to \( {\sigma} \in [-\tfrac{1}{2}, 0] \times [0, \infty) \), since \( B(\e, \o, \s_1, \s_2)\) is conjugate to \(  B(\e, \o, -\s_1, \s_2)\) via the symmetry operators \( R_j , j=1,2, \) defined by 
 \begin{equation*} 
(R_1W)(x)=W(-x) \quad \text{and} \quad (R_2W)(x)=\overline{W}(x).
 \end{equation*}
In both cases we have 
\begin{equation*}
B(\varepsilon, \omega, \sigma_1, \sigma_2) = R_j B(\varepsilon, \omega, -\sigma_1, \sigma_2) R_j,  \quad j=1, 2. 
\end{equation*}
Consequently, it  suffices to examine the eigenvalue problems of $B(\e, \o, \s)$ on $L^2(\mathbb T)$:
\be \label{EEB_r}
0=[B(\e, \o, \sigma)-\l]W
\ee
for $\s \in [-\frac{1}{2}, 0] \times [0, \infty)$, rather than for all $\s \in \RR^2$. 

\begin{figure}
    \centering
    \includegraphics[width=0.56 \textwidth]{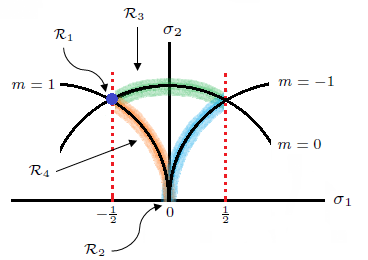}
    \caption{The four regions $\mathcal{R}_1$ -- $\mathcal{R}_4$. The blue dot indicates $\sigma^*=(-\frac{1}{2}, \frac{\sqrt{3}}{2})$. Figure reproduced from \cite{CJ}.} 
    \label{unstable wave vectors}
\end{figure}

To characterize unstable Bloch wave vectors $\s \in [-\frac{1}{2}, 0] \times [0, \infty)$, we first analyzed the eigenvalue problem \eqref{EEB_r} at the critical value $\e=0$. For small $\e>0$ the operators $B(\eps,\omega, \s)$ can be regarded as small perturbations of
\be
B(0,\omega,\s)=-(1+(\partial_x+i\s_1)^2-\s_2^2)^2.
\notag
\ee
Since this operator has constant coefficients, its eigenvalues $\mu_m(\s)$ are explicitly given by
\be \label{eigenvalues}
\mu_m(\s)=-(1-(m+\sigma_1)^2-\sigma_2^2)^2 \leq 0, \quad m \in \ZZ.  
\ee
Therefore, for small $\e>0$ unstable Bloch wave vectors $\s$ can occur when 
\be 
(m+\sigma_1)^2+\sigma_2^2 \approx1, 
\notag
\ee
from which we readily identify the set 
\be \label{S0}
\mathcal{S}_0:=\{ \s \in [-\frac{1}{2}, 0] \times [0, \infty)   ~|~ C_0: \s_1^2+\s_2^2=1, ~ C_{\pm 1}: (\s_1\pm1)^2+\s_2^2=1\}.  
\ee
We then distinguish four regions around the set $\mathcal{S}_0$: for some $\delta>0$, 
\be
\begin{split}
\mathcal{R}_1: & = \{ \s \in [-\frac{1}{2}, 0] \times [0, \infty) ~| ~dist(\s, \sigma^*) \leq \delta \},  \\
\mathcal{R}_2: & = \{ \s \in [-\frac{1}{2}, 0] \times [0, \infty) ~| ~dist(\s, (0,0)) \leq \delta \},\\
\mathcal{R}_3: & =  \{ \s \in [-\frac{1}{2}, 0] \times [0, \infty) ~| ~dist(\s, \mathcal{C}_0) \leq \delta,  ~ dist(\s, \sigma^* ) \geq \delta \},  \\
\text{and} \quad \mathcal{R}_4: & =  \{ \s \in [-\frac{1}{2}, 0] \times [0, \infty) ~| ~dist(\s, \mathcal{C}_1) \leq \delta,  ~dist(\s, \sigma^*) \geq \delta, ~   dist(\s, (0,0)) \geq \delta \}.
\notag
\end{split}
\ee
Here $\sigma^*=(-\frac{1}{2}, \frac{\sqrt{3}}{2})$. Notice that the two circles $C_0$ and $C_1$ (or $C_1$ and $C_{-1}$, respectively) intersect at $\sigma^*=(-\frac{1}{2}, \frac{\sqrt{3}}{2})$ (or $(0,0)$, respectively); see Figure \ref{unstable wave vectors}.  In each of the regions $\mathcal{R}_1$ and $\mathcal{R}_2$, the operator $B(0, \o, \s)$ has two small eigenvalues, whereas in $\mathcal{R}_3$ and $\mathcal{R}_4$ it has only one.  In \cite{CJ}, by analyzing the eigenvalue problem \eqref{EEB_r} in $\mathcal{R}_1$, we established Theorem \ref{spectral instability}.

\smallskip

On the other hand, if $\sigma$ is bounded away from $\mathcal{S}_0$, the eigenvalues $\mu_m(\s)$ are uniformly negative. Hence the eigenvalues of $B(\eps,\omega,\sigma)$, being small perturbations of those of $B(0,\omega,\sigma)$, also remain uniformly negative  provided $\eps$ is sufficiently small. In our previous work \cite{CJ}, it was enough to identify the spectrally unstable Bloch wave vectors, and the spectral behavior away from $\mathcal{S}_0$ was mentioned without detailed analysis. In the present work, however, such localization is essential for establishing the existence of the most unstable eigenvalue. We therefore state the precise result as a lemma below, with the proof deferred to Appendix $A$.
\begin{lemma} \label{sqrte}
 Let $s$, $b$ and $\tilde \e_0$ be chosen as in Theorem \ref{spectral instability}. Then there exists $0 <\eta = O(\sqrt{\e})$  such that for all $\e \in (0, {\tilde \e}_0]$ and all $\o \in (-\frac{1}{2}, \frac{1}{2})$, 
\be
\mathcal S_{\e, \o} \subset \{ \s : \textrm{dist}(\s,  S_0) < \eta \}. 
\notag
\ee
That is, all unstable Bloch wave vectors are contained in a band of width $O(\sqrt{\e})$ around $S_0$. 
\end{lemma}

\medskip

We now concern the existence of a maximal eigenvalue and the local analyticity of the  associated spectral branch near the maximizer.  In the lemma below, we use term  \textit{a domain of analyticity}:  $\calR \subset [-1/2, 1/2]\times \bbr$ is a domain of analyticity if each point in $\calR$ has a neighborhood where $\la$ is represented as an analytic function.

\begin{lemma}\label{analytic lambda_3}
 Let $s$, $b$ and $\tilde \e_0$ be chosen as in Theorem \ref{spectral instability}. There exists  ${\bar \e}_0 \in (0, \tilde \e_0]$ such that for each $0<\e \leq \bar \e_0$ and each $\o \in (\frac{1}{2}, \frac{1}{2})$, the Bloch operator $B(\e, \o, \sigma)$ admits the most unstable eigenvalue $\l_M$, attained at some $\bar \sigma \in \RR^2$. Moreover, there exists a neighborhood $\mathcal{N} \subset \RR^2$ of $\bar \sigma$ and a real analytic function 
\be
\l: \mathcal{N}  \longrightarrow \RR
\notag
\ee
such that $\l(\bar \sigma)=\l_M$ and for each $\sigma \in \mathcal{N} $, $\l(\sigma)$ is an unstable eigenvalue of $B(\e, \o, \sigma)$. In particular,  $\l(\sigma)$ defines a locally analytic branch of unstable eigenvalues selected near $\bar\sigma$.
\end{lemma}

\begin{proof} 
Fix any $\o \in (\frac{1}{2}, \frac{1}{2})$. For each $\s_0 \in \mathcal{S}_0$, an unstable eigenvalue $\lambda(\e, \o, \sigma)$ of $B(\e, \o, \sigma)$ for $(\e, \sigma)$ near $(0, \sigma_0)$ was constructed in \cite{CJ} via the Lyapunov–Schmidt reduction, after separating the four regions $\mathcal{R}_1-\mathcal{R}_4$. More precisely, $B(\e, \o, \sigma)$ has at most one or two simple unstable eigenvalues, since $B(0, \o, \s)$ has two critical eigenvalues for $\s \in \mathcal{R}_1$ and $\s \in \mathcal{R}_2$, and only one for $\s \in \mathcal{R}_3$ and $\s \in \mathcal{R}_4$. 

The analyticity of $\lambda$ in $(\e, \sigma)$ near $(0, \sigma_0)$ follows from two ingredients of that analysis: the implicit function theorem applied in the Lyapunov–Schmidt procedure, and the Weierstrass preparation theorem characterizing the leading-order behavior of $\l$  in $(\e, \sigma)$ near $(0, \sigma_0)$  (see, for instance, the second equation in $(47)$ and $(60)$  in \cite[Section 3.2]{CJ} for $\s \in \mathcal{R}_1$). Since $S_0$ is compact, there exists $\gamma > 0$ such that  for all $\e \in (0, \tilde \e_0)$ the domain of analyticity of  $\lambda$ contains 
\be
\calR^{\gamma}:= \{ \sig\in [-\frac 12, 0] \times \bbr ~|~ dist(\sig, \mathcal S_0) < \gamma\}.
\ee 
Here we note that $\gamma$ depends only on the fixed constant $\e_0$ from Theorem \ref{existence_roll}. Accordingly, the Lyapunov–Schmidt reduction yields at most one or two analytic branches of unstable eigenvalues. However, the choice between these branches is unambiguous, since one branch is always greater than or equal to the other.

By Lemma \ref{sqrte}, all unstable Bloch wave vectors are contained in a band of width $O(\sqrt{\e})$ around $S_0$. We now choose  $\bar \e_{0} \in (0, \tilde \e_0]$ sufficiently small so that 
\be
O (\sqrt {\bar \e_0}) < \gamma. 
\notag
\ee
Then, for any $\e \in (0, \bar \e_0]$, the corresponding unstable Bloch wave vectors $\sigma$ lie entirely within $\calR^{\gamma}$. Consequently, every unstable eigenvalue, including the largest one, is contained in an analytic neighborhood of $\lambda$ centered at some $\sigma_0$ on $S_0$, and therefore admits a locally analytic representation in $\sigma$. 
\end{proof}

\smallskip

We now construct an exponentially growing solution to the linearized equation by replacing the unstable eigenvalue $\lambda(\sigma^*)$ in \cite[Proposition 3.1]{CJ} with the largest unstable eigenvalue $\lambda_M = \lambda(\bar\sigma)$ obtained in Lemma \ref{analytic lambda_3}. This replacement not only captures the dominant growth behavior but also provides an upper bound estimate on the growing solution, which will play a key role in controlling the error terms in the subsequent nonlinear analysis.

To proceed, let $s$, $b$ and $\bar \e_0$ be chosen as in Lemma \ref{analytic lambda_3}, and fix any $\e \in (0, \bar \e_0]$ and $\omega \in (-\frac 12, \frac 12)$. For notational convenience, we suppress the explicit dependence of $\mathcal{L}_{\e, \o}$ and $B(\e, \o, \s)$, and simply write $\mathcal{L}$ and $B(\s)$, respectively. In addition, we define a small square of Bloch wave vectors centered around the maximizer $\bar\sigma = (\bar\sigma_1, \bar\sigma_2)$:
\begin{equation}\label{ItimesJ}
I \times J := [\bar\sigma_1, \bar\sigma_1 + c_0] \times [\bar\sigma_2, \bar\sigma_2 + c_0],
\end{equation}
where $c_0 > 0$ is chosen sufficiently small so that the square $I \times J$ lies entirely within the analyticity region $U$ (a quantitative bound to be given later in Section \ref{Nonlinear instability}). This guarantees that the unstable eigenvalue $\lambda(\sigma)$ remains real--analytic throughout $I \times J$.


\begin{proposition}\label{liin}
Let $W(x, \sig)$ be the continuous family of the eigenfunctions of $B(\sigma)$ corresponding to the unstable eigenvalue $\la(\sig)$ for $\sig \in I \times J$. Define the function 
\be \label{egsolution}
U(x, y) = \int_{I\times J} e^{i (\sig_1 x + \sig_2 y)} W(x, \sig) d\sig_1 d\sig_2.
\ee
Then the semigroup evolution $e^{\cL t } U$ satisfies the estimate
\be
\| (e^{\cL t } U)(\cdot) \|_{H^2(\RR^2)} \simeq e^{\l_Mt} e(t),
\notag
\ee
where the function $e(t)$ satisfies
\be \label{upperlower}
\frac{C_1}{(1+ t)^{\frac 1l}} \le  e(t) \le C_2 
\ee
for some positive constants $C_1$, $C_2$ and $l$, depending on $\eps, \omega$, $\bar \sig$ and $c_0$. 
\end{proposition}
\begin{proof}
We define a function 
\be \label{soperatorV}
V(t, x, y): = \int_{I\times J} e^{\la(\sig) t} e^{i(\sig_1 x+\sig_2 y)} W(x, \sig) d\sig_1 d\sig_2. 
\ee
Evaluating at $t=0$ gives $V(0, x, y) = U(x, y)$. Since $\l$ is bounded on the compact set $I \times J$ and $W(x, \cdot) \in L^1(I \times J)$, direct differentiation under the integral sign yields
\be
\partial_t V = \int_{I\times J} \la(\sig) e^{\la(\sig) t}  e^{i(\sig_1 x+\sig_2 y)} W(x, \sig) d\sig_1 d\sig_2. 
\notag
\ee
On the other hand, using the relation \eqref{EEB_r}, we compute
\be
\begin{split}
\cL V 
& = \int _{I\times J}e^{\la(\sig) t}  e^{i(\sig_1 x+\sig_2 y)} B(\sigma)W(x, \sig) d\sig_1 d\sig_2 \\
& =\int_{I\times J} \la(\sig) e^{\la(\s) t}  e^{i(\sig_1 x+\sig_2 y)}W(x, \sig) d\sig_1 d\sig_2,
\notag
\end{split}
\ee
which confirms that $V$ indeed satisfies the linearized equation with initial data $U$, i.e., $e^{\cL t }U = V$. 

To estimate the norm of $V$, we apply Plancherel's Theorem (see Lemma \ref{norm_W} in the end of the section), yielding 
\be\label{L2grow}
\|V(t, \cdot)\|^2_{H^2(\RR^2)}
\simeq \int_J \int_I e^{2\la(\sig)t} \left  \|  W(\cdot, \s) \right\|_{H^2(\mathbb T)}^2 d\s_1 d\s_2. 
\ee
By the continuity of the mapping $\sig \in I\times J \mapsto W(\cdot, \sig) \in H^2(\mathbb T)$, there exist constants  $C_1$ and $C_2$, depending on $\e$, $\omega$ and $\bar \sigma$, such that 
\be \label{V_H2}
{ C_1 \int_J \int_I  e^{2\la(\sig)t}  d\sig_1 d\sig_2  \le \|V(t, \cdot)\|^2_{H^2(\bbR^2)} \le C_2 \int_J \int_I  e^{2\la(\sig)t} d\sig_1 d\sig_2. }
\ee
 Since  $\la(\sig)$ is analytic on $I\times J$, there exists $0 < l < \infty$ such that the following expansion about $\s=\bar \s$ holds:
\be\label{diff}
 \la(\sig)-\la(\bar \sigma) =\sum_{k=0}^l a_{k}(\sig_1-\bar \sigma_1)^{l-k}(\sig_2 - \bar \sigma_2)^{k} + o(|\sig-\bar \sigma|^l),
 \ee
where not all $a_k$ vanish. Substituting this expansion into the integrand and letting $\sig_i- \bar \sigma_i = \kappa_i$ for $i=1, 2$,  yields that
\begin{align*}
\int_J \int_I  e^{2\la(\sig)t}  d \sig_1 d \sig_2 
  &= 
   e^{2\la_M t}  (e(t))^2,
   \end{align*}
   where we define
   \be \label{et}
    e(t):=  \Big(\int ^{c_0}_0\int^{c_0}_0 
  e^{\left(2\sum_{k=0}^l a_{k}{\kappa_1}^{l-k}{\kappa_2}^{k}+ o(|\kappa|^l)\right)t} d \kappa_1 d\kappa_2 \Big)^{1/2}. 
  \ee
Recalling $\l_M=\l(\bar \sig)$ is the maximal growth rate, the right--hand side of \eqref{diff} is nonpositive from which we deduce that
\be \label{et_ul}
\Big(\int ^{c_0}_0 \int^{c_0}_0 e^{-C (\kappa_1^l + \kappa_2^l)t} d\kappa_1 d\kappa_2 \Big)^{1/2}
\le   e(t) \le C
\ee
 for some constant $C>0$.

We now estimate the lower bound in \eqref{et_ul} by $c_0e^{-Cc_0^l}$ for $ 0\le t\le 1$, and by
$C(c_0) t^{ -\frac 1l}$ for $ t>1$. In particular, for the case $t>1$ we perform the change of variables $p_i = \kappa_i^l  t$ for $i=1,2$, which gives
\be
\int_0^{ c_0} e^{-C\kappa_1^lt} d\kappa_1 = t^{-\frac 1l} \int_0^{c_0^lt} \frac{p_1^{\frac 1l -1}e^{-Cp_1}}{l} dp_1 
\notag
\ee
Putting these estimates together, we arrive at 
 \[ C_1\frac{e^{\la_M t}}{(1+ t)^{\frac{1}{l}}} \le \|V(t, \cdot)\|_{L^2(\bbR^2)} \simeq e^{\la_M t} e(t) \le C_2e^{\la_M t}\]
 for some positive constants $C_1$, $C_2$ depending on  $\e$, $\omega$, $\bar \sigma$ and $c_0$. 
\end{proof}

\begin{remark} \label{remark_e(t)}
It is worth emphasizing the essential difference between the present two-dimensional setting and the one-dimensional analysis in \cite{JLL}. In the one-dimensional case, the expansion \eqref{diff} near the maximally unstable mode admits uniform estimate of the form
\be\label{1dcase}
-2C(\s-\bar \s)^l  \leq \l(\s)-\l(\bar \s) \leq -\frac{1}{2}C(\s-\bar \s)^l, 
\ee
which lead to precise decay:
\be
\| (e^{\cL t } U)(\cdot) \|_{H^2(\RR)} \simeq \frac{e^{\l_Mt}}{(1+t)^{1/l}}. 
\notag
\ee
In contrast, in the two-dimensional problem considered here, the bound \eqref{1dcase} generally fail due to the more complicated geometry near its maximum. As a result, $e(t)$ defined in \eqref{et} only satisfies the bounds \eqref{upperlower}. For this reason, throughout the nonlinear analysis in Section \ref{Nonlinear instability}, we consistently express all estimates in terms of the combined factor $e^{\l_Mt} e(t)$, as it represents the natural growth scale imposed by the linear dynamics. 

In addition, we note that in Proposition \ref{liin}, the function $e(t)$ is decreasing in $t$, and satisfies
\be \label{de_e}
|e'(t)| \leq Cc_0^le(t)
\ee
for some constant $C >0$. In what follows $Cc_0^l$ can be made arbitrarily small by taking $c_0$ sufficiently small. This fact will play an important role in controlling several constants in the preceding sections. For example, we will choose $c_0$ such that 
\be \label{C(c0)3}
\l_M-Cc_0^l>0
\ee
from which we deduce that the function  $e^{\l_Mt}e(t)$ is increasing in $t$, because
\be
\frac{d}{dt}\Big[e^{\l_Mt}e(t)\Big]=\l_M e^{\l_Mt}e(t)+e^{\l_Mt}e'(t) \geq (\l_M-Cc_0^l)e^{\l_Mt}e(t)>0.
\notag
\ee
This monotonicity will also play a key role in the nonlinear instability estimates in Section \ref{Nonlinear instability}. 
\end{remark}

In the end of the section, we introduce a lemma stating the Parseval's theorem holds for the Bloch transformation,
extending \cite[Lemma 3.3]{JLL} to the two--dimensional case. 
\begin{lemma}\label{norm_W}  For nonnegative integer $p$ and any $W(x, \sigma) \in H^p(\mathbb T)$,
\be
\Big\| \int_{I\times J} e^{i (\s_1 x + \s_2 y)} W(x, \s) d\s \Big\|^2_{H^p(\RR^2)} \simeq \int_{I\times J} \|W(\cdot, \s)\|^2_{H^p(\mathbb T)} d \s.
\notag
\ee
Here, the Sobolev space $H^p(\mathbb T)$ with periodic boundary conditions can be defined as 
\be  \label{norm}
H^p(\mathbb T)=\{ h \in L^2(\mathbb T) ~|~ \|h\|_{H^p(\mathbb T)}^2:=\sum_j (1+j^2)^p |\hat h(j)|^2 < \infty  \},
\ee
where $\hat h(j):=\frac{1}{2\pi} \int_0^{2\pi} e^{-ijx} h(x) dx$, $j \in \ZZ$, referred to as the Fourier coefficient of $h$. In particular, $\|h\|^2_{L^2(\mathbb T)}:=\sum_j |\hat h(j)|^2$. 
\end{lemma}
\begin{proof}
For each $j \in \ZZ$, we denote
\be
I_j \times J:=   [\bar\sigma_1+j , \bar\sigma_1 + c_0 +j] \times
[\bar\sigma_2, \bar\sigma_2+ c_0]. 
\notag
\ee
Since $W(x, \s)$ is $2\pi$--periodic in $x$, $W$ has a Fourier expansion 
\be \label{fs}
W(x, \s)=\sum_{j \in \ZZ} \hat W(j, \s) e^{ijx}. 
\ee   
Denoting by $k=(k_1, k_2)$ a multiindex of order $0 \leq |k|=k_1+k_2 \leq p$, it follows that  
\be \label{eW}
\begin{split}
& \partial^{k}\int_{I\times J} e^{i (\s_1 x + \s_2 y)} W(x, \s) d\s \\
& =  \int_{I\times J} \sum_{j \in \ZZ} (i\s_2)^{k_2}  (i(\s_1+j))^{k_1}  e^{i((\s_1+j)x+\s_2 y)}\hat W(j, \s) d\s \\
& =\sum_{j \in \ZZ} \int_{I_j \times J}  (i\s_2)^{k_2}  (i\s_1)^{k_1}  e^{i(\s_1x+\s_2 y)}\hat W(j, \s_1-j, \s_2) d\s_1 d\s_2 \\
& = \int_{\RR^2}  e^{i(\s_1x+ \s_2 y)} \sum_{j \in \ZZ} (i\s_2)^{k_2}  (i\s_1)^{k_1} \chi_{I_j \times J}(\s) \hat W(j, \s_1-j, \s_2) d\s_1 d\s_2.
\end{split}
\ee
Here, $ \chi_{I_j \times J}$ is the characteristic function on $I_j \times J$. Notice that since $c_0<<1$, $\chi_{I_j \times J}\cdot \chi_{I_k \times J}=0$ unless $j=k$, and the right--hand side of \eqref{eW} is an inverse Fourier transform of the function $\sum_{j \in \ZZ}(i\s_2)^{k_2}  (i(\s_1+j))^{k_1}  \chi_{I_j \times J}(\s) \hat W(j, \sig_1-j, \sig_2)$. Thus we deduce from Parseval's identity and \eqref{norm} that  for $0 \leq |k| \leq p$, 
\be
\begin{split}
&  \|  \partial^{k} \int_{I\times J} e^{i (\s_1 x + \s_2 y)} W(x, \s) d\s \|^2_{L^2((x, y); \RR^2)} \\
& \simeq  \| \sum_{j \in \ZZ} (i\s_2)^{k_2}  (i\s_1)^{k_1}   \chi_{I_j \times J}(\s) \hat W(j, \s_1-j, \s_2)  \|^2_{L^2(\s; \RR^2)} \\
& =  \int_{\RR^2} \sum_{j \in \ZZ} |\s_2|^{2k_2} |\s_1|^{2k_1} | \chi_{I_j \times J}(\s)  |^2   | \hat W(j, \s_1-j, \s_2)|^2 d\s_1 d\s_2 \\
& = \sum_{j \in \ZZ}  \int_{I_j \times J} |\s_2|^{2k_2} |\s_1|^{2k_1} | \hat W(j, \s_1-j, \s_2)|^2 d\s_1 d\s_2 \\
& = \int_{I \times J}  \sum_{j \in \ZZ}  |\s_2|^{2k_2} |\s_1+j|^{2k_1} | \hat W(j, \s) |^2 d\s\\
& \simeq \int_{I \times J}  \sum_{j \in \ZZ}(1+j^2)^{k_1} | \hat W(j, \s) |^2 d\s.
\end{split}
\notag
\ee
Noting $0 \leq k_1 \leq p$, it completes the proof. 
\end{proof}

\section{Semigroup estimates}\label{Semigroup estimates}

In this section, we obtain semigroup estimates for the linearized evolution operator $e^{t\mathcal{L}}$, where  $\mathcal{L}$ is defined in \eqref{linearization about rolls}. 
 In particular, we demonstrate \eqref{asdg_L} that the semigroup norm is governed by the maximal growth rate of unstable eigenvalues of $\calL$, which provide us  a growth bound for solutions to an inhomogeneous equation (see Proposition \ref{gr_of_inhomo}).

To analyze $e^{t\mathcal{L}}$, we briefly review the Bloch transform framework (see, e.g., \cite{JNRZ1, JmZk, JZ}).
For the general theory on the Bloch transform, we refer to
\cite[Section XIII.16]{RS4}.  For any $w \in L^2(\RR^2)$, the inverse Fourier transform formula yields the representation 
\be \label{semi1}
w(x, y) =\frac{1}{2\pi}\int_{-\infty}^{\infty} \int_{-\frac 12}^{\frac 12} e^{i( \s_1x + \s_2 y)} \check w (x, \s)  d \s_1 d\s_2.
\ee
We define  $\check w $, the Bloch transform of $w$,  in $L^2(\T, [-\frac 12, \frac 12]\times \bbr)$ by
 $$\check w (x, \s): = \sum_{j \in \ZZ} e^{ijx} \hat v (\s_1+j, \s_2),$$ where $\hat w$ denotes the $2$--dimensional Fourier transform of $w$. 
 From the Parseval's identity $L^2$ isometry holds,
 \[ \| w\|_{L^2(\bbr^2)}^2 = \int_{-\infty}^{\infty}\int_{-\frac{1}{2}} ^{\frac{1}{2}} \| \check w (\cdot, \s)\|^2_{L^2(\T)} d\s_1 d\s_2.\]
 Moreover the $H^p$ isometry is proven in Lemma \ref{norm_W},
 \be\label{norm_W_id}
\| w \|_{H^p(\RR^2)}^2 \simeq \int_{-\infty}^{\infty} \int_{-\frac 12}^{\frac 12} \| \check w (\cdot, \s) \|^2_{H^p(\mathbb T)} d\s_1 d\s_2.
\ee
We can write the linear evolution by  $e^{t\mathcal L}$ in terms of the associated Bloch operators defined in  \eqref{Bloch_r}: 
\be \label{semi2}
e^{t\mathcal L}w(x, y) = \int_{-\infty}^{\infty} \int_{-\frac 12}^{\frac 12} e^{i( \s_1x + \s_2 y)} e^{tB(\s)} \check w (x, \s)  d \s_1 d\s_2.
\ee
 It together with \eqref{norm_W_id} leads
 \be \label{norm_eq}
\| e^{t\mathcal L}w\|^2_{H^p(\RR^2)} \simeq \int_{-\infty}^{\infty} \int_{-\frac 12}^{\frac 12} \| e^{tB(\s)}\check w (\cdot, \s)\|^2_{H^p(\mathbb T)} d\s_1 d\s_2.
\ee
Therefore, to estimate $e^{t\mathcal{L}}$, it suffices to study the semigroup $e^{tB(\sigma)}$. 
The proof follows a similar line of reasoning as in the classical result that the heat semigroup 
$e^{t\Delta}$ is compact on $H_0^2(\Omega)$ for a bounded domain $\Omega$ 
(see \cite[Section~22.4]{Lax}). 
\smallskip


\begin{lemma}\label{semigroup_compact}
For any fixed $\s \in [-\frac 12, \frac 12] \times \bbr$ and time $T>0$ the semigroup $$e^{TB(\s)}: {H^2(\mathbb T)} \to {H^2(\mathbb T)}$$ is a compact operator.  
\end{lemma}
\begin{proof} 
Recalling the Bloch operator \eqref{Bloch_r}, for notational simplicity we set
\be\label{pqB_0}
D_x:= \pa_{x} + i\sigma_1, \quad  \mathcal{A}_0:=(1+2\varepsilon\omega)D_x^2 +1-\si_2^2, \quad \CalF(\tilde u):= \e^2 +2b\tilde u_{\e, \o}- 3s \tilde u^2_{\e, \o}. 
\ee
Then the solution $u:=e^{tB(\sigma)} w$ satisfies
\begin{equation}\label{eqB0}
 \pa_t u = -\mathcal{A}_0 u + \CalF(\tilde u) u, \qquad u(0,\cdot)=w. 
 \end{equation}
Here we note that $D_x$ plays the role of a covariant-type derivative, and integration by parts using $D_x$ follows rules analogous to those for ordinary derivatives:  for $f, g \in C^{\infty}(\mathbb T)$ and $\a, \b \in \CC$,
\begin{align*}
\int_{\mathbb T} \bar f D_x g dx & = - \int_{\mathbb T} \overline{D_x f}g dx,\\
\int_{\mathbb T} \bar f D^2_x g dx & = - \int_{\mathbb T} \overline{D_x f}D_x g dx,\\
\int_{\mathbb T} \bar f (\a D^2_x  + \b)g  dx & =  \int_{\mathbb T} \overline{(\a D_x^2 + \b)  f} g dx.
\end{align*}

Multiplying \eqref{eqB0} by $\bar u$, integrating resultant equality in $x$ over $\mathbb T$ and applying integration by parts, we obtain 
\begin{equation}\label{l2_eq}
\frac{1}{2}\ddt \|u(t)\|_{L^2(\mathbb T)}^2 =   -\|\mathcal{A}_0 u(t)\|_{L^2(\mathbb T)}^2 +\langle \mathcal{F}(\tilde u)u, u \rangle_{L^2(\mathbb T)} 
\end{equation}
We fix a finite time $T>0$. Since $|\CalF(\tilde u)|\leq C$, applying the Gronwall inequality yields that for all $t \in [0, T]$, 
\begin{equation}\label{l2_bdd}
\| u(t)\|^2_{L^2(\mathbb T)} \le C(T)\| w\|^2_{L^2(\mathbb T)}
\end{equation}
for some $C(T)$ depending on $T$. 

Next, multiply \eqref{eqB0} by $t \overline{\mathcal{A}_0^2u}$ and integrate the resultant equality in $x$ over $\mathbb T$. Then integration by parts gives
\[ \frac t2  \frac{d}{dt} \|\mathcal{A}_0u(t)\|_{L^2(\mathbb T)}^2 = - t\|\mathcal{A}_0^2u(t)\|_{L^2(\mathbb T)}^2 
+   \Re \langle t\mathcal{F}(\tilde u)u, \mathcal{A}_0^2u \rangle_{L^2(\mathbb T)} \]
We now integrate the above equation in $t$ over $[0, T]$. Apply integration by parts to the left-hand side, and use  the Young's inequality and \eqref{l2_bdd} to the right-hand side to obtain that  
\begin{align*}
& \frac 12 \left [T \|\mathcal{A}_0u(T)\|_{L^2(\mathbb T)}^2 - \int_0^T \|\mathcal{A}_0u(t)\|_{L^2(\mathbb T)}^2  dt \right] \\
\qquad & = -   \int_0^T \|\mathcal{A}_0^2u(t)\|_{L^2(\mathbb T)}^2 dt
+  \Re \int_0^T  \langle t\mathcal{F}(\tilde u)u, \mathcal{A}_0^2u \rangle_{L^2(\mathbb T)} dt\\
\qquad &  \le   C(T) \| w\|^2_{L^2(\mathbb T)}. 
\end{align*}
Due to \eqref{l2_eq} and \eqref{l2_bdd}, it follows that
\begin{align*}
\frac T2 \|\mathcal{A}_0u(T)\|_{L^2(\mathbb T)}^2 & \le \frac 12 \int_0^T \|\mathcal{A}_0u(t)\|_{L^2(\mathbb T)}^2 dt  + 
C(T) \| w\|^2_{L^2(\mathbb T)} \\
 & \leq  -\frac 14\int_0^T \ddt \|u(t)\|_{L^2(\mathbb T)}^2 dt+  C(T) \| w\|^2_{L^2(\mathbb T)}  \\
 & \leq C(T) \| w\|^2_{L^2(\mathbb T)}.
\end{align*}
Consequently, recalling \eqref{pqB_0}, we see that 
\be \label{D2}
\|D^2_xu(T)\|^2_{L^2(\mathbb T)}   \le C(T, \s_2)\| w\|^2_{L^2(\mathbb T)}.
\ee
Here, by direct calculation
\begin{align*}
|D^2_xu|^2=|\pa^2_x u|^2  + 4\sig_1^2|\pa_x u |^2 + \sig_1^4 |u|^2 -2\sig_1^2 \Re(\pa_x^2u \bar u) + 4\sig_1^3 \Im(\pa_x u \bar u) + 4\sig_1 \Im(\pa_x^2 u \overline{\pa_x u})
\end{align*}
Applying the Young's inequality  and using \eqref{l2_bdd} and \eqref{D2}, we have
\be\label{secondDu}
 \|\pa_x^2 u(T)\|^2_{L^2(\mathbb T)}
 \le C\| \pa_x u(T)\|^2_{L^2(\mathbb T)}+C(T, \s_2)\| w\|^2_{L^2(\mathbb T)}. 
\ee 
By plugging this inequality into the Gagliardo-Nirenberg inequality:
 \[ \| \pa_x u(t)\|_{L^2(\mathbb T)} \le C \|u(t)\|^{\frac 12}_{L^2(\mathbb T)} \| \pa_x^2 u(t)\|^{\frac 12}_{L^2(\mathbb T)},\]
we deduce that 
 \[  \| \pa_x u(T)\|^4_{L^2(\mathbb T)} \le C(T, \s_2)  \|w\|^2_{L^2(\mathbb T)}(\|\pa_x u(T)\|^{2}_{L^2(\mathbb T)}+\|w\|^2_{L^2(\mathbb T)}).\]
Simplifying further, we arrive at 
\[ \| \pa_x u(T)\|_{L^2(\mathbb T)}^2  \le C(T, \s_2)\| w\|^2_{L^2(\mathbb T)}.\]
Combining this estimate with \eqref{l2_bdd} and \eqref{secondDu} we conclude that
\[ \|u (T)\|_{H^2(\mathbb T)}^2 \le C(T, \s_2)\| w\|^2_{L^2(\mathbb T)}.\]
Higher-order estimates can be derived in a similar manner. The compactness of the semigroup $e^{TB(\sigma)}$ on $H^2(\mathbb T)$ follows from the above estimates and the compact Sobolev embedding.
\end{proof}







\smallskip

We now demonstrate that the semigroup generated by the linearized operator exhibits exponential growth, with rate bounded by the largest unstable spectrum. 

\begin{lemma}\label{semigroup bound} Let $\l_M=\l(\bar \sigma)$ be the maximum growth rate defined in Lemma \ref{analytic lambda_3}. Then there exists a constant $C$ 
such that  for any $w \in H^2(\RR^2)$, 
\be  \label{asdg_L}
\| e^{t \mathcal L}w\|_{H^2(\RR^2)} \leq Ce^{\l_M t} \|w\|_{H^2(\RR^2)}. 
\ee

\end{lemma}
\begin{proof} Let $t>0$. We first notice that $e^{tB(\s)}: {H^2(\mathbb T)} \to {H^2(\mathbb T)}$ is the strongly continuous $1$--parameter semigroup. Moreover, by the previous lemma, $e^{tB(\s)}$ is compact. Since  $B(\s)$ is self-adjoint in $L^2(\T)$, so is $e^{tB(\s)}$ (see \cite[Theorem 4.1.1] {A}). Then, by the spectral theorem for compact self-adjoint operators and the spectral mapping theorem (\cite[Theorem 13. Section 34.5]{Lax}), the spectrum of $e^{tB(\s)}$ consists only of eigenvalues of the form $e^{\l_j t}$, where $\l_j \in spec(B(\s))$. Writing $W=\sum \a_j W_j$ in the orthonormal eigenbasis $\{W_j \}$ of $L^2(\T)$, we obtain 
\be
\| e^{tB(\s)}W \|_{L^2(\mathbb T)}^2=\| \sum \a_j e^{\l_j t}W_j \|_{L^2(\mathbb T)}^2 \leq e^{2\l_M t} \|W \|_ {L^2(\mathbb T)}^2.
\notag
\ee
The estimate \eqref{asdg_L} in $L^2(\RR^2)$ follows directly from the Bloch decomposition and Parseval's identity (see \eqref{norm_W_id} and \eqref{norm_eq}). For the higher estimate, we set 
\be \label{mathcalA}
\mathcal A:=1+(1+2\e \o)\partial_x^2+\partial_y^2. 
\ee
Then the solution $u:=e^{t \calL }w$ satisfies
\be \label{23}
\pa_t u =-\mathcal A^2u+ \CalF(\tilde u)u, \qquad u(0,\cdot)=w,
\ee
where $\CalF(\tilde u)$ is defined in \eqref{pqB_0}. 

Noting that $H^2(\RR^2)$-norm is equivalent to the combination of $\|u \|_{L^2(\RR^2)}$ and $\|Au \|_{L^2(\RR^2)}$, it suffices to estimate $\|Au \|_{L^2(\RR^2)}$. We multiply \eqref{23} by $\mathcal A^2 u$ and integrate over $\bbr^2$. Using the Young's inequality and the bound $|\CalF(\tilde u)|\leq C$, we obtain
\be
\begin{split}
\ddt \| Au(t)\|_{L^2(\bbr^2)}^2 
& =-2 \| A^2u (t)\|^2_{L^2(\bbr^2)} +2\langle  \CalF(\tilde u)u, A^2u \rangle_{L^2(\bbr^2)} \\
& \leq C \|u(t)\|_{L^2(\bbr^2)}^2 \\
& \leq C e^{2\l_Mt} \|w\|_{L^2(\bbr^2)}^2.
\end{split}
\notag
\ee
Here we used the estimate \eqref{asdg_L} in $L^2(\mathbb T)$ for the last inequality. Integrating the above relation with respect to time $t$ yields 
\be
\begin{split}
\| Au(t)\|_{L^2(\bbr^2)}^2  
& \leq \| Au(0)\|_{L^2(\bbr^2)}^2 + C(\|w\|_{L^2(\bbr^2)}^2 + e^{2\l_M t}\|w \|_{L^2(\bbr^2)}^2) \\
& \leq  C(\|w \|_{H^2(\bbr^2)}^2 + e^{2\l_M t}\|w \|_{L^2(\bbr^2)}^2) \\
& \leq C e^{2\l_M t}\|w\|_{L^2(\bbr^2)}^2, 
\end{split}
\notag
\ee
which completes the proof.
\end{proof}

\medskip

As a key application of the semigroup estimate, we now derive a growth bound for solutions to an inhomogeneous linear equation. This result will later be used to control higher-order nonlinear terms in the instability analysis.

\begin{proposition} \label{gr_of_inhomo}Let $e(t)$ be as defined by \eqref{et} in Proposition \ref{liin}. Suppose the inhomogeneous forcing term $g(t,x,y)$  satisfies
\be \label{nonhomo}
\| g(t)\|_{H^2(\RR^2)} \leq K \big(e^{\l_M t} e(t)\big)^j, \quad t \geq 0, \quad 2 \leq j \leq N
\ee
for some $K>0$. Then the solution $u(t,x,y)$ to the initial value problem
\be \label{equnonhomo}
 \qquad \partial_t u = \mathcal L u + g, \quad u(0,\cdot)=0
\ee
satisfies the estimate
\be \label{gr}
\qquad \|u(t)\|_{H^2(\RR^2)} \leq CK \big(e^{\l_M t} e(t)\big)^j, \quad t \geq 0
\ee 
for some constant $C>0$. 
\end{proposition}
\begin{proof} By the Duhamel formula and the estimate \eqref{asdg_L}, together with the assumption \eqref{equnonhomo}, we obtain
\be \label{H2norm}
\begin{split}
 \| u(t)\|_{H^2(\RR^2) } \leq \int_0^t \| e^{(t-s)\mathcal{L}} g(s) \|_{H^2(\RR^2)} ds
 & \leq CK e^{\l_M t} \int_0^t  e^{(j-1)\l_M s} e(s)^j ds.
\end{split}
\ee
Recalling the constant $Cc_0^l$ defined in \eqref{de_e}, we apply integration by parts to obtain
\be
\begin{split}
 \int_0^t  e^{(j-1)\l_M s} e(s)^j ds 
 \leq \frac{1}{(j-1)\l_M} \left(e^{(j-1)\l_M t} e(t)^j +Cjc_0^l  \int_0^te^{(j-1)\l_M s} e(s)^{j}ds \right). 
\notag
\end{split}
\ee
Rewriting, we have 
\be
\Big(1-\frac{Cjc_0^l}{(j-1)\l_M} \Big)   \int_0^t  e^{(j-1)\l_M s} e(s)^j ds \leq \frac{1}{(j-1)\l_M} e^{(j-1)\l_M t} e(t)^j.
\notag
\ee
We now choose $c_0>0$ small enough such that 
\be \label{C(c0)1}
0<\frac{Cjc_0^l}{(j-1)\l_M}  \leq \frac{2Cc_0^l}{\l_M}<\frac{1}{2},
\ee
which yields
\be
\int_0^t  e^{(j-1)\l_M s} e(s)^j ds \leq  \frac{2}{(j-1)\l_M} e^{(j-1)\l_M t} e(t)^j \leq \frac{2}{\l_M} e^{(j-1)\l_M t} e(t)^j. 
\notag
\ee
Substituting this into the right-hand side of \eqref{H2norm} gives
\be 
 \| u(t)\|_{H^2(\RR^2) } \leq \frac{C K}{\la_M}   (e^{\l_Mt}  e(t))^j. 
\notag
\ee
\end{proof}


\section{Nonlinear instability}\label{Nonlinear instability}

In this section, we prove the main theorem of the paper, i.e., the nonlinear instability of the roll solutions $\tilde u(x)$ to the two-dimensional generalized Swift-Hohenberg equation \eqref{2DgSH_w}. More precisely,  for a fixed small $\delta>0$, we define an appropriate positive constant $\theta$ and construct a solution $V_\delta (t, x, y)$ to \eqref{2DgSH_w} that satisfy \eqref{goal_delta} and \eqref{goal_theta}. 

Following the strategy developed in \cite{Gr, JLL, RT1}, we decompose the perturbed solution into two parts: an approximate solution that captures the dominant linear growth and a remainder term that accounts for the nonlinear corrections. More precisely, we consider the ansatz, 
\be \label{V_delta}
V_{\delta}(t, x, y) = V^{\textup{app}}(t, x,y) + v(t, x, y), 
\ee
where $V^{\textup{app}}(t,x,y)$ is constructed to approximate the evolution of the unstable mode up to a prescribed order in $\delta$, and $v$ denotes the nonlinear remainder which is expected to remain small over the relevant time scale. 

The approximate solution takes the form
\be \label{Vapp}
V^{\textup{app}}(t, x, y) = \tilde u(x) + \sum_{j=1}^N\delta^{j}V_j(t, x, y). 
\ee
The definitions of $V_j$ and $v$ are provided in \eqref{ex} and \eqref{veq}, respectively.
The first-order term $V_1$ is the exponentially growing solution of the linear problem with the maximal growth rate, and higher-order corrections $V_j$ for $j\geq 2$ are designed to compensate for  the nonlinear self-interaction of $V_1$. These terms are solutions of recursively defined inhomogeneous linear equations derived from lower-order contributions. Their inclusion ensures that the approximate solution matches the true dynamics up to order $\mathcal{O}(\delta^N)$, where the integer $N$ will be chosen subsequently.

The initial data corresponding to \eqref{V_delta} is taken as 
\be\label{initial perturbation}
V_{\delta}(0, x, y)=\tilde u(x)+\delta U(x, y),
\ee
where $U(x,y)$ denotes the initial data associated with the unstable linear mode constructed in Proposition \ref{liin}, that is, $V_1(0, x,y)=U(x,y)$. 
Besides $V_1$ we impose 
\be
V_j(0, x, y)=0 \quad \text{for all $j\geq 1$}, \quad \text{and} \quad v(0,x,y)=0. 
\notag
\ee

In what follows, recalling \eqref{mathcalA} and \eqref{pqB_0}, we rewrite the original equation \eqref{2DgSH_w} and its linearized operator \eqref{linearization about rolls}  in the form
\be\label{2DgSH_w_r}
\partial_t u=Lu+f(u) \quad \text{and} \quad \mathcal L=L+\mathcal F(\tilde u),
\ee
where $L:=-\mathcal A^2$ and $f(u):=\e^2 u+bu^2-su^3$.

An $N$-th order approximate solution $V^{\textup{app}}(t, x, y)$ is constructed as follows. We first set $$V_1:=e^{t\mathcal L}U$$ as defined in Proposition \ref{liin}. That is, $V_1=V_1(t, x, y)$ satisfies 
\be \label{V_1}
\partial_t V_1-\mathcal LV_1=0, \quad V_1(0, \cdot)=U. 
\ee
To define the higher-order terms $V_j$ for $2 \leq j \leq N$,  we insert the ansatz \eqref{Vapp} and compute 
\be \label{eq_V}
\begin{split}
& \partial_t V^{\textup{app}} -LV^{\textup{app}} - f(V^{\textup{app}}) \\
 &=  \partial_t V^{\textup{app}} -LV^{\textup{app}} - f(\tilde u) - (f(V^{\textup{app}}) - f(\tilde u))\\
&= \partial_t \left( \sum_{j=2}^N\delta^{j}V_j \right) - \mathcal L \left( \sum_{j=2}^N\delta^{j}V_j \right) +(3s \tilde u-b) \left( \sum_{j=1}^N\delta^{j}V_j \right)^2 + s \left( \sum_{j=1}^N\delta^{j}V_j \right)^3. 
\end{split}
\ee
Here we used $L\tilde u+f(\tilde u)=0$ and \eqref{V_1}, and expanded $f(V^{\textup{app}})$ around $\tilde u$: 
\be
f(V^{\textup{app}})-f(\tilde u)= f'(\tu)(\va-\tu) +\frac 12 f''(\tu)(\va-\tu)^2 + \frac {1}{3!}f'''(\tu)(\va-\tu)^3,
\notag
\ee
\be
 f'(\tu)=\mathcal{F}(\tu)=\e^2+2b\tu-3s\tu^2, \quad f''(\tu)= 2b-6s\tu, \quad f'''(\tu)= -6s.
\notag
\ee
The nonlinear terms on the right-hand side of \eqref{eq_V} can be organized by power of $\delta$
\begin{align*}
\left( \sum_{j=1}^N\delta^{j}V_j \right)^2 & =\sum_{j=2}^N\delta^j  \left(\sum_{k+p=j}  V_kV_p \right)+\sum_{k+p=N+1}^{2N} \delta^{k+p} V_k V_p,  \\
\left( \sum_{j=1}^N\delta^{j}V_j \right)^3 & =\sum_{j=3}^N \delta^j \left(\sum_{k+p+q=j} V_kV_pV_q\right)+\sum_{k+p+q=N+1}^{3N} \delta^{k+p+q} V_k V_p V_q
\end{align*}
with  $1 \leq k, p, q \leq N$. Grouping terms of order $\delta^j$, we define $V_j$ for $2 \leq j \leq N$ recursively by the linear inhomogeneous problems
\begin{align} \label{ex} \begin{aligned}
\partial_t V_2 - \mathcal L V_2  & = -(3s \tilde u-b) V_1^2, \\
\partial_t V_j - \mathcal L V_j &= -(3s \tilde u-b) \sum_{k+p=j}  V_kV_p -s\sum_{k+p+q=j} V_kV_pV_q \quad \mbox{for }3\le j\le N
\end{aligned}
\end{align}
with  $V_j(0, \cdot)=0$ for $2\le j\le N$. Substituting these equations back into \eqref{eq_V}, the approximate solution satisfies
\be\label{mathcal_R}
\begin{split}
 \partial_t V^{\textup{app}} -LV^{\textup{app}} - f(V^{\textup{app}}) & = (3s \tilde u-b) \sum_{k+p=N+1}^{2N} \delta^{k+p} V_k V_p + s \sum_{k+p+q=N+1}^{3N} \delta^{k+p+q} V_k V_p V_q
 \\
 & =: \mathcal{R}^{\textup{app}}
\end{split}
\ee
with $1\le k, p, q \le N$.
The symbol $ \mathcal{R}^{\textup{app}}$ represents the residual of the approximate solution. 

Finally, inserting the ansatz \eqref{V_delta} into \eqref{2DgSH_w_r} yields the equation for the remainder $v$:
\begin{align}\label{veq} \begin{aligned}
&\partial_t v 
= Lv +\big(\e^2+2b \va-3s(\va)^2 \big)v+ (b-3s\va )v^2-sv^3+ \mathcal{R}^{\textup{app}},\\
&v(0, \cdot) = 0.
\end{aligned}\end{align}

\medskip

The profile of $j$-th correction term $V_j$ can be estimated as follows.

\begin{lemma}\label{estimateV}
There exists a uniform constant $\rho >0$  such that
\be \label{estV}
\|V_j(t)\|_{H^2(\RR^2)} \leq (\rho e^{ \l_M t}e(t)\big)^j , \quad   1 \leq j\leq N. 
\ee
\end{lemma}
\begin{proof}
We first claim that  there exists a positive sequence $\{C_j\}_{1\le j\le N}$ satisfying
\[\|V_j(t)\|_{H^2(\RR^2)} \leq C_j(e^{ \l_M t}e(t)\big)^j , \quad   1 \leq j\leq N .\]
By  Proposition \ref{liin} we have 
\be \label{estimateV1}
\|V_1(t) \|_{H^2(\RR^2)} \leq C_1 e^{\l_M t} e(t)
\ee 
for some $C_1>0$, so the claim holds for $j=1$. Assume now that the claim holds for all $V_l$ with $1\le l \le j-1$ for some $2\le j\le N$. By the induction hypothesis,
\begin{align*}
\| \sum_{k+p=j}  V_kV_p(t)\|_{H^2(\bbr^2)} \le \sum_{k+p=j}  \|V_k (t)\|_{H^2(\bbr^2)} \|V_p (t)\|_{H^2(\bbr^2)} \le \sum_{k+p =j} C_kC_p  (e^{ \l_M t}e(t)\big)^j,
\end{align*}
and similarly, 
\begin{align*}
\|\sum_{k+p+q=j} V_kV_pV_q (t) \|_{H^2(\bbr^2)} &\le \sum_{k+p+q=j}  \|V_k(t) \|_{H^2(\bbr^2)} \|V_p (t)\|_{H^2(\bbr^2)} \|V_q (t)\|_{H^2(\bbr^2)}\\
&\le 
\sum_{k+p+q=j} C_kC_pC_q(e^{ \l_M t}e(t)\big)^j
\end{align*}
with $1\le k, p, q \le N$. For notational simplicity, we set
\be
\a:=C\|3s\tilde u -b \|_{L^\infty(\RR)}, \quad \tilde s:=Cs
\notag
\ee
where $C$ is defined in \eqref{gr}.  Applying Proposition \ref{gr_of_inhomo} to \eqref{ex} yields
\[\|V_j(t)\|_{H^2(\RR^2)} \leq C_j (e^{ \l_M t}e(t)\big)^j,\]
where the sequence $\{C_j\}_{2\le j\le N}$ is defined recursively by 
\begin{align*}
 C_2 & :=  \a C_1^2, \\
 C_j&:= \a \sum_{k+p =j} C_kC_p + \tilde s \sum_{k+p+q=j} C_kC_pC_q,  \quad j\ge 3.
\end{align*}
In fact, there exists $\rho>0$ such that
\be\label{catalan}
C_j \le \rho^j.
\ee  
When $\tilde s=0$, the recurrence relation reduces to that of the Catalan numbers, which grow exponentially. The proof of \eqref{catalan} is postponed to Appendix~B.
\end{proof}

\medskip

Due to the above lemma, the residual term $ \mathcal{R}^{\textup{app}}$ in \eqref{mathcal_R}  can be estimated by
\be
\|\mathcal{R}^{\textup{app}} (t)\|_{H^2(\RR^2)} \leq CN^2 \big(\delta \rho e^{\l_Mt}e(t)\big)^{N+1}(1+\delta \rho  e^{\l_M t}e(t)+ \cdots).
\notag
\ee
Here, $N^2$ is the upper bound for $|\{(k, p, q)| k+p+q=N, 1\le  k, p, q \le N-1\}|$. 

Now fix $\theta \in (0, \frac{1}{2})$ to be specified later. Recalling that $e^{\l_Mt}e(t)$ is increasing in $t$ for sufficiently small $c_0$ satisfying \eqref{C(c0)3}, we define $T^{\delta}>0$ as the finite time such that 
\be \label{Time}
\delta \rho e^{\l_M T^{\delta}}e(T^{\delta}) =\theta.  
\ee
Then, for all $0\le t\le T^{\delta}$, the residual term $\mathcal{R}^{\textup{app}}$ satisfies
\begin{align}\label{estimate_R}
\|  \mathcal{R}^{\textup{app}} (t)\|_{H^2(\RR^2)} 
 \le CN^2\frac{ \big(\delta \rho e^{\l_Mt}e(t)\big)^{N+1}}{1-\theta} 
 \le CN^2\big(\delta \rho e^{\l_Mt}e(t)\big)^{N+1}. 
\end{align}
For later use, we record the following consistency condition on $c_0$:
\be \label{C(c0)4}
\delta \rho e(0)=\delta \rho c_0^2  \leq \theta. 
\ee

\medskip

The next lemma required to prove Theorem \ref{nonlinear instability} concerns the energy estimates for the remainder $v$ up to time $T^{\delta}$.
The local well-posedness of $v$ in $H^2(\bbr^2)$ follows from a standard fixed point argument and is therefore omitted.  We define the maximal time $T_M$ by
\be\label{TM}
T_M:=\sup \Big\{T>0  ~|~  \|v (t) \|_{H^2(\RR^2)} \leq \frac{1}{2} \quad \text{for all $0 \leq t \leq T$} \Big\}. 
\ee
In the case where the bound holds for all times, we set $T_M=+\infty$.  Note that the maximal existence time is larger than $T_M$.

\begin{lemma} [Energy estimate] \label{Energy estimate} ~
\begin{enumerate}
\item  For $T^{\delta}$ defined in \eqref{Time} it holds that $T^\delta < T_M$.
\item  The remainder $v$ satisfies the same bound as $\mathcal{R}^{\textup{app}}$; that is, for all $ 0\le t\le T^{\delta}$,
\be \label{Rv} 
\|v (t) \|_{H^2(\RR^2)} \leq  CN^2\big(\delta \rho  e^{\l_Mt}e(t)\big)^{N+1}.
\ee
\end{enumerate}
\end{lemma}
\begin{proof}
To prove  point $(1)$ by contradiction, suppose that $T_M \leq T^{\delta}$.  Recalling \eqref{mathcalA} and \eqref{2DgSH_w_r}, 
the $v$–equation \eqref{veq} can be rewritten as
\[ \pa_t v = -\mathcal{A}^2v +P(v) +\rap,\]
where $ P(v) =(\e^2+2b \va-3s(\va)^2)v+ (b-3s\va)v^2-sv^3$.  Then integration by parts yields
\be
\begin{split}
\ddt \| v (t)\|_{L^2(\bbr^2)}^2 & =-2\| Av (t)\|^2_{L^2(\bbr^2)} + 2\langle P(v) +\rap, v \rangle_{L^2(\bbr^2)}, \\
\ddt \| Av (t)\|_{L^2(\bbr^2)}^2 & =-2 \| A^2v (t)\|^2_{L^2(\bbr^2)} +2\langle  P(v) +\rap, A^2v \rangle_{L^2(\bbr^2)}.
\end{split}
\notag
\ee
Recalling \eqref{Vapp}, \eqref{estV} and \eqref{TM},  we have that for all $0< t \le T_M$, 
\be
\begin{split}
\| v(t)\|_{L^{\infty}(\bbr^2)} & \le C, \\
\|  \va(t) \|_{L^{\infty}(\bbr^2)} & \le C \big(\e +\sum_{j=1}^N (\delta \rho  e^{\l_M t} e(t))^j \big) \le C(\e + \theta). 
\end{split}
\notag
\ee
It then follows that for all $0< t \le T_M$, 
\[ \| P(v)(t)\|_{L^2(\bbr^2)} \le C\|v(t)\|_{L^2(\bbr^2)}. \]
Consequently, 
\be
\ddt \left(\| v(t)\|^2_{L^{2}(\bbr^2)} +  \| Av (t)\|^2_{L^2(\bbr^2)}  \right) 
\le C \left( \| v(t)\|^2_{L^{2}(\bbr^2)} + \| \rap(t) \|^2 _{L^2(\bbr^2)} \right). 
\notag
\ee
Applying Gronwall inequality and using $v(0,\cdot)=0$ together with \eqref{estimate_R}, we obtain for all $0 \leq t \leq T_M$,
\begin{align} \label{Gronwall} \begin{aligned}
\|v (t) \|^2_{H^2(\RR^2)} 
& \leq \int_0^t e^{C(t-s)}\| \mathcal{R}^{\textup{app}}(s) \|^2_{H^2(\RR^2)}ds \\
& \leq CN^4 (\delta \rho)^{2(N+1)}e^{Ct}\int_0^t  e^{\big(2(N+1)\l_M -C\big)s}e(s)^{2(N+1)} ds. 
\end{aligned}
\end{align}
We now choose $N$ sufficiently large so that 
\be \label{chooseN}
 2(N+1)\l_M -C >0.
\ee
As in the proof of Proposition \ref{gr_of_inhomo}, the integral on the right-hand side of \eqref{Gronwall} can be calculated by 
\be \label{inte}
\int_0^t  e^{\big(2(N+1)\l_M -C\big)s}e(s)^{2(N+1)} ds \leq \frac{2}{2(N+1)\l_M -C}  e^{\big(2(N+1)\l_M -C\big)t}e(t)^{2(N+1)}
\ee
provided that $c_0$ is chosen sufficiently small so that 
\be \label{C(c0)2}
\frac{2(N+1)Cc_0^l}{2(N+1)\l_M -C}<\frac{1}{2}. 
\ee 
Substituting \eqref{inte} inito \eqref{Gronwall} gives, for all $0 \leq t \leq T_M$, 
\be \label{v_bound} 
\|v (t) \|_{H^2(\RR^2)} \leq CN^2 \big(\delta \rho e^{\l_Mt}e(t)\big)^{N+1} \leq CN^2\theta^{N+1}. 
\ee
By continuity, the maximal time for which $\|v (t) \|_{H^2(\RR^2)} \leq \frac{1}{2}$ can be extended beyond $T_M$ if $\theta$ is chosen such that 
\be \label{smalltheta}
CN^2\theta^{N+1} <\frac{1}{2}. 
\ee
This contradicts the definition of $T_M$, and thus completes the proof of $(1)$. As a consequence, by repeating the same argument up to time $T^{\delta}$, the estimate \eqref{v_bound} immediately yields point $(2)$.
\end{proof}

Equipped with the previous lemmas, we are now in a position to prove the nonlinear instability stated in Theorem \ref{nonlinear instability}.

\begin{proof}[Proof of Theorem \ref{nonlinear instability}] 
From \eqref{estV} and \eqref{Rv}, the deviation of the full solution from $\tilde u$ at time $T^{\delta}$ can be estimated as
\be \label{deviation}
\begin{split}
\| V_{\delta}(T^{\delta}) - \tilde u\|_{L^2(\RR^2)}  & = \| \sum_{j=1}^N \delta^j V_j (T^{\delta})+ v (T^{\delta})\|_{L^2(\RR^2)} \\
& \geq  \| \delta V_1 (T^{\delta}) \|_{L^2(\RR^2)} - \| \sum_{j=2}^N \delta^j V_j (T^{\delta})\|_{L^2(\RR^2)} - \|v (T^{\delta})\|_{L^2(\RR^2)} \\
&  \geq C\theta-  2C\theta^2 - CN^2 \theta^{N+1}\\
& \ge \frac{C}{2}\theta,
\end{split}
\ee
where the last inequality holds for $\theta \in (0, \frac{1}{2})$ sufficiently small, still satisfying the smallness condition \eqref{smalltheta}. Lastly, based on 
\eqref{C(c0)3}, \eqref{C(c0)1}, \eqref{C(c0)4} and \eqref{C(c0)2}, we choose $c_0>0$ small enough to meet all the preceding requirements.

We now verify that $T^\delta=O(|\ln \delta|)$. From the defining relation \eqref{Time} and the bounds in  \eqref{upperlower}, we have
\be \label{3}
C_1 \frac{\delta \rho e^{\l_Mt} }{(1+t)^\frac 1l} \le \delta \rho e^{\l_Mt}e(t) \le C_2\delta \rho e^{\l_M t}.
\ee
Let $t_1$ and $t_2$ denote the unique times satisfying
\be
C_1 \frac{\delta \rho e^{\l_Mt_1} }{(1+t_1)^\frac 1l} =\theta \quad \text{and} \quad C_2\delta \rho e^{\l_M t_2}=\theta.
\notag
\ee
Since all three functions in \eqref{3} are increasing in $t$, it follows that  
$$t_2 \leq T^{\delta} \leq t_1.$$
Both $t_1$ and $t_2$ clearly satisfy $O(|\ln \delta|)$, and thus we conclude that $$T^{\delta}=O(|\ln \delta|).$$
\end{proof}

\begin{remark} Throughout the analysis, the constants $N$, $\theta$, and $c_0$ were chosen in the following order:
first, $N$ was fixed large enough to satisfy \eqref{chooseN}; next, $\theta \in (0,\tfrac12)$ was selected sufficiently small to meet the smallness condition \eqref{smalltheta} and the last inequality of \eqref{deviation}; finally, $c_0>0$ was chosen small enough to ensure \eqref{C(c0)3}, \eqref{C(c0)1}, \eqref{C(c0)4} and \eqref{C(c0)2}. This hierarchical choice guarantees that all preceding inequalities are compatible.
\end{remark}

\smallskip

\section{Appendix A}\label{Appendix A} 

 In this appendix we prove Lemma \ref{sqrte}. The argument repeatedly uses computations from the proof of Lemma \ref{semigroup_compact}, in particular integration by parts for covariant derivatives.

\smallskip

\begin{proof}[Proof of Lemma \ref{sqrte}] 

Let $\e \in (0, \tilde \e_0]$ and $\o \in (-\frac{1}{2}, \frac{1}{2})$, where $\tilde \e_0$  is chosen as in Theorem \ref{spectral instability}. For a $\eta \in (0, 1]$, we define 
\[\cG^{\eta} := \{\s \in [-\frac 12, 0] \times [0, \infty)~ | ~\textrm{dist}(\sigma, S_0) \ge \eta\}.\] 
From the defining formula \eqref{eigenvalues}, since
the set $\{ \mu_m(\sigma) | m\in \mathbb Z\}$ is discrete,  we have that for $\s \in  \cG^{\eta}$,
\[  \mu_m(\s) \leq  
\begin{cases}  
-\frac{25}{16} ~~ \mbox{ for } |m| \geq 2,  \\
-\eta^2~~ \mbox{ for } m=0, \pm 1,
\end{cases}\]
therefore  $\mu_m(\sigma)$ has a negative upper bound $-\eta^2$.

We first estimate the operator norm $\| B(\e, \o, \s)- B(0, \o, \s)\|_{L(L^2(\mathbb T), L^2(\mathbb T))}$. Proceeding as in the proof of Lemma \ref{semigroup_compact}, we have \begin{align}\label{Op}\begin{aligned}
& \langle B(\e, \o, \s)u- B(0, \o, \s)u, u \rangle_{L^2(\mathbb T)} 
\\
& \leq  4\e|\o| \left( |1-\s_2^2|\| D_x u\|_{L^2(\mathbb T)}^2 +\| D^2_x u\|_{L^2(\mathbb T)}^2\right) + \langle \mathcal{F}(\tilde u)u, u \rangle_{L^2(\mathbb T)}
\end{aligned}\end{align}
 for $u \in H^2(\mathbb T)$.

 Let $V \in H^2(\mathbb T)$ and $\lam \in \bbr$ solve the eigenvalue equation 
\begin{equation*}
B(\e, \o, \s)V =\la V
\end{equation*}
 for $\s \in \cG^{\eta}$.
Multiplying $\bar V$ to the above equation and integrating by parts, we have
\be\label{identy_appen}
\| [(1 +2\e\o)D_x^2  + (1-\s_2^2) ]V\|^2_{L^2(\mathbb T)} = \langle (\mathcal{F}(\tilde u)-\l)V,  V \rangle_{L^2(\mathbb T)}.
\ee

\medskip

We now fix $M >2$ and consider the following two cases to show that $\l <0$. 

\smallskip

\noindent Case I. $0\le \s_2 \le M$ \\
Substituting the Young's inequality
$$2\Re \langle (1+2\e\o)D_x^2V,  (1-2\s_2^2)V \rangle_{L^2(\mathbb T)} \geq - \frac 14 (1 +2\e\o)^2 \| D_x^2 V\|^2 _{L^2(\mathbb T)}- 4 ( 1-\s_2^2 )^2 \| V\|^2_{L^2(\mathbb T)}$$
into \eqref{identy_appen} and recalling $\cF(\tilde u) =O(\e)$, we have
\begin{equation*}
   (1 +2\e\o )^2 \| D_x^2 V\|^2_{L^2(\mathbb T)} \le 4 ( 1-\s_2^2 )^2 \|V\|^2_{L^2(\mathbb T)} + \frac{4}{3}(O(\e) -\lam)\|V\|^2_{L^2(\mathbb T)}.
\end{equation*}
It implies that
\be 
\| D^2_x V \|^2_{L^2(\mathbb T)} \lesssim (M^4 + (\e -\lam))\|V\|^2_{L^2(\T)}.
\notag
\ee
Consequently we have

\begin{align*}
\lambda \| V \|^2_{L^2(\mathbb{T})}
&= \langle B(\varepsilon, \omega, \sigma) V, V \rangle_{L^2(\mathbb{T})} \\
&= \langle (B(\varepsilon, \omega, \sigma)- B(0, \omega, \sigma))V, V \rangle_{L^2(\mathbb{T})}
+ \langle B(0, \omega, \sigma)V, V \rangle_{L^2(\mathbb{T})} \\
&\lesssim \varepsilon |\omega| (M^4 + (\varepsilon - \lambda)) \| V \|^2_{L^2(\mathbb{T})}
+ \varepsilon \| V \|^2_{L^2(\mathbb{T})}
-\eta^2 \| V \|^2_{L^2(\mathbb{T})}.
\end{align*} 
Here we used the interpolation  $\|D_xV\|^2_{L^2(\T)} \lesssim \|V\|_{L^2(\T)}\|D^2_x V\|_{L^2(\T)}$ together with $\mu_m(\s) \leq -\eta^2$. Finally, choosing $\eta^2= O( \e)$ such that 
$ \e|\o|(M^4 + \e)+ \e < \eta^2$, we ensure that  $\la <0$.

\medskip

\noindent Case II. $ \s_2 \ge M$  \\
 We note that
the cross term on the left-hand side of \eqref{identy_appen} is positive:
 \[2(1 +2\e\o) (1-\s_2^2) \Re \langle D_x^2V, V \rangle_{L^2(\mathbb T)}  = -2(1 +2\e\o) (1-\s_2^2)\| D_x V\|^2_{L^2(\mathbb T)} >0.\]
Plugging this into \eqref{identy_appen} gives
 \[ (1+ 2\e\o)^2 \| D_x^2 V\|^2_{L^2(\T)} +2(1 +2\e\o) (\s_2^2-1)\| D_x V\|^2_{L^2(\mathbb T)}  +(1-\s_2^2)^2\|V\|^2_{L^2(\T)} \leq ( O(\e) - \lam)\|V\|^2_{L^2(\T)}, \]
hence, if we choose $\tilde \e_0 < \frac 12$ \textit{e.g.}, we bound the left-hand side of \eqref{Op} by 
 \[(16\e |\o|(O(\e)-\lambda) + \e)\| V\|_{L^2(\mathbb T)}^2.\]
Similarly as in Case I, we obtain
\begin{align*}
\lambda \| V \|^2_{L^2(\mathbb{T})} 
& =  \langle (B(\varepsilon, \omega, \sigma) - B(0, \omega, \sigma))V, V \rangle_{L^2(\mathbb{T})} + \langle B(0, \omega, \sigma)V, V \rangle_{L^2(\mathbb{T})} \\
&\lesssim  (16\e |\o|(O(\e)-\lambda) + \e)\| V\|_{L^2(\mathbb T)}^2- \eta^2 \| V \|^2_{L^2(\mathbb{T})}.
\end{align*}
Choosing $\eta^2= O( \e)$ such that $ 16 \e O(\e)|\o|+ \e < \eta^2$, we ensure that  $\la <0$.

\end{proof}

\smallskip

\section{Appendix B}
In this section we provide the proof of \eqref{catalan} used in proving Lemma \ref{estV}.
This argument is a simple application of the generating-function method widely used in enumerative combinatorics (see, e.g., \cite{FO1990} and \cite{FS}). 
\begin{proposition}
For a given $C_1>0$ define the sequence $(C_j)_{j\ge2}$ by
\begin{align*}
C_2&=\alpha\,C_1^2,\\
C_j&=\alpha\!\!\sum_{k+p=j} C_k C_p\;+\;\tilde s\!\!\sum_{k+p+q=j} C_k C_p C_q,\qquad j\ge3,
\end{align*}
where $k, p, q \ge1$. Then there exists a constant $\rho>0$ such that
\begin{equation}\label{eq:expbd}
C_j \le \rho^{\,j}\qquad\text{for all }j\ge1.
\end{equation}
\end{proposition}

\begin{proof}
Let $C(z):=\sum_{j\ge1} C_j z^j$. The recurrence relations imply
\[
C(z) \;=\; C_1 z \;+\; \alpha\,C^2(z) \;+\; \tilde s\,C^3(z)^ .
\]
Equivalently, with 
\[
F(C,z):= C - C_1 z - \alpha C^2 - \tilde s C^3,
\]
we have $F\bigl(C(z),z\bigr)\equiv0$. Note that $F(0,0)=0$ and
\[
\partial_C F(0,0)=1\neq0.
\]
By the implicit function theorem, there is $R>0$ such that $C(z)$ is analytic on the disk $\{|z|<R\}$; in particular $C$ has a (strictly) positive radius of convergence $R$.

Since $C_1>0$ and $\alpha,\tilde s\ge0$, the recursion shows inductively that all $C_j\ge0$. Hence, by the Cauchy--Hadamard formula,
\[
\limsup_{j\to\infty} C_j^{1/j}\;\le\; R^{-1}.
\]
Therefore, for any $\varepsilon>0$ there exists $j_0$ such that $C_j\le (R^{-1}+\varepsilon)^j$ for all $j\ge j_0$. Increasing a base we can reach a constant $\rho>0$ with \eqref{eq:expbd} holding for all $j\ge1$.
\end{proof}
\begin{remark}
When $\tilde s=0$, the recurrence relations  defines the Catalan sequence. In this case $C_j \sim \frac{ (4\alpha)^{j-1}C_1^j}{\sqrt{\pi} j^{\frac 32}}$  is well known. 
\end{remark}


\end{document}